\documentclass[11pt,openbib]{article}

\usepackage{amsmath}
\usepackage{amssymb}
\usepackage{amsfonts}
\usepackage[T1]{fontenc}
\usepackage{graphicx}

\newcommand{\C}{\mathbb{C}}

\newcommand{\Z}{\mathbb{Z}}

\newcommand{\dee}{\mathrm{d}}
\newcommand{\comp}{\, \raisebox{2pt}{$\scriptstyle\circ \, $}}
\newcommand{\rowspace}{\rule{0pt}{16pt}}

\newcommand{\lefthook}{\mbox{$\, \rule{8pt}{.5pt}\rule{.5pt}{6pt}\, \, $}}

\newcommand{\setrule}{\, \rule[-4pt]{.5pt}{13pt}\, }

\allowdisplaybreaks

\begin{document}

\title{\textbf{Shifting operators in geometric quantization}}
\author{Richard Cushman and J\k{e}drzej \'{S}niatycki \thanks{%
Department of Mathematics and Statistics, University of Calgary. \newline
email: rcushman@ucalgary.ca and sniatycki@ucalgary.ca }}
\date{}
\maketitle

\addtocounter{footnote}{1} \footnotetext{%
printed: \today}

\section{Introduction}

%%%%%%%%%%%%%%%

In a series of papers on Bohr-Sommerfeld-Heisenberg quantization of
completely integrable systems \cite{cushman-sniatycki1}, \cite%
{cushman-sniatycki2}, \cite{cushman-sniatycki3}, \cite{sniatycki15}, we
introduced operators acting on eigenstates of the action operators by
shifting the eigenvalue of the $j^{\mathrm{th}}$ momentum by $\hbar $ and
multiplying the wave function by $e^{-i\theta _{j}}$. Here $(I_{j},\theta
_{j})$ are classical action angle coordinates; that is, the angle $\theta
_{j}$ is a multivalued function determined up to $2n\pi $. The need for such
operators in the Bohr-Sommerfeld quantization was pointed out already by
Heiseberg \cite{heisenberg25}. However, there was no agreement about their
existence and their interpretation. The aim of this paper is to derive
shifting operators from the first principles of geometric quantization.
\medskip

In the following, we use rescaled action angle coordinates $(j_{i},\vartheta
_{i})$, where the angle $\vartheta _{j}$ is a multivalued function
determined up to an arbitrary integer $n$, defined on an open domain $U$ of $%
P$. In these coordinates, the symplectic form of our phase space $(P,\omega
) $ has local expression ${\omega }_{\mid U} =\sum_{i=1}^{k}{\mathrm{d} j}_i
\wedge \mathrm{d} {\vartheta}_i$. \medskip

The first step of geometric quantization is called prequantization. It
consists of the construction of a complex line bundle $\pi :L\rightarrow P$
with connection whose curvature form satisfies a prequantization condition
relating it to the symplectic form $\omega $. A comprehensive study of
prequantization, from the point of view of representation theory, was given
by Kostant in \cite{kostant}. The work of Souriau \cite{souriau} was aimed
at quantization of physical systems, and studied a circle bundle over phase
space. In Souriau's work, the prequantization condition explicitly involved
Planck's constant $h$. In 1973, Blattner \cite{blattner73} combined the
approaches of Kostant and Souriau by using the complex line bundle with the
prequantization condition involving Planck's constant. Since then, geometric
quantization has been an effective tool in quantum theory. \medskip

We find it convenient to deal with connection and curvature of line bundles
in the framework the theory of principal and associated bundles 
\cite{kobayashi-nomizu}. In this framework, the prequantization condition reads 
\begin{equation*}
\mathrm{d} \beta =(\pi ^{\times })^{\ast }(-%
\mbox{${\scriptstyle
\frac{{1}}{{h}}}$} \, \omega ),
\end{equation*}%
where $\beta $ is the connection $1$-form on the principal 
$\mathbb{C}^{\times }$-bundle $\pi ^{\times }:L^{\times }\rightarrow P$ associated to
the complex line bundle $\pi :L\rightarrow P$, and $\mathbb{C}^{\times }$ is
the multiplicative group of nonzero complex numbers. \medskip

The aim of prequantization is to construct a representation of the Poisson
algebra $(C^{\infty }(P), \{ \, \, , \, \, \} , \cdot )$ of $(P,\omega )$ on the space of sections of the line bundle $L$. Each Hamiltonian vector field $X_{f}$ on $P$ lifts to a
unique $\mathbb{C}^{\times }$-invariant vector field $Z_{f}$ on $L^{\times }$
that preserves the principal connection $\beta $ on $L^{\times }$. If the
vector field $X_{f}$ is complete, then it generates a 1-parameter group $%
\mathrm{e}^{tX_{f}}$ of symplectomorphisms of $(P,\omega )$. Moreover, the
vector field $Z_{f}$ is complete and it generates a $1$-parameter group $%
\mathrm{e}^{tZ_{f}}$ of connection preserving diffeomorphisms of the bundle $%
(L^{\times },\beta )$, called quantomorphisms, which cover the $1$-parameter
group ${\mathrm{e}}^{tX_{f}}$.\footnote{%
The term quantomorphism was introduced by Souriau \cite{souriau} in the
context of $\mathrm{SU}(n)$-principal bundles and discussed in detail in his
book \cite{souriau70}. The construction discussed here follows \cite%
{sniatycki80}, where the term quantomorphism was not used.} In this case, ${%
\mathrm{e}}^{tX_{f}}$ and ${\mathrm{e}}^{tZ_{f}}$ are 1-parameter groups of
diffeomorphisms of $P$ and $L^{\times }$, respectively. We shall refer to ${%
\mathrm{e}}^{tX_{f}}$ and ${\mathrm{e}}^{tZ_{f}}$ as flows of $X_{f}$ and $%
Z_{f}$. Since $L$ is an associated bundle of $L^{\times }$, the action ${%
\mathrm{e}}^{\,tZ_{f}}:L^{\times }\rightarrow L^{\times }$, induces an
action ${\widehat{\mathrm{e}}}^{\,tZ_{f}}:L\rightarrow L,$ which gives rise
to an action on smooth sections $\sigma $ of $L$ by push forwards, $\sigma
\mapsto {\widehat{\mathrm{e}}}_{\ast }^{\,tZ_{f}}\sigma ={\widehat{\mathrm{e}%
}}^{\,tZ_{f}}\,\raisebox{2pt}{$\scriptstyle\circ \, $}\sigma \,%
\raisebox{2pt}{$\scriptstyle\circ \, $}{\mathrm{e}}^{\,-tX_{f}\text{.}}$.
Though ${\widehat{\mathrm{e}}}_{\ast }^{\,tZ_{f}}\sigma $ may not be defined
for all $\sigma $ and all $t$, its derivative at $t=0$ is defined for all
smooth sections. The prequantization operator 
\begin{equation}
\mathcal{P}_{f}\sigma =i\hbar \, 
\mbox{${\displaystyle \frac{\dee }{\dee t}}
\rule[-10pt]{.5pt}{25pt} \raisebox{-10pt}{$\, {\scriptstyle t=0}$}$} \mathrm{%
\widehat{e}}_{\ast }^{\, tZ_{f}}\sigma ,  \label{1.1}
\end{equation}%
where $\hbar $ is Planck's constant divided by $2\pi $, is a symmetric
operator on the Hilbert space $\mathfrak{H}$ of square integrable sections
of $L$. The operator $\mathcal{P}_{f}$ is self adjoint if $X_{f}$ is
complete.\medskip

The whole analysis of prequantization is concerned with \emph{globally}
Hamiltonian vector fields. Since every vector field on $(P,\omega )$ that
preserves the symplectic form is locally Hamiltonian, it is of interest to
understand how much of prequantization can be extended to this case. In
particular, we are interested in is the case where the Hamiltonian vector fields are 
the vector fields of the angle variables $\vartheta _{i}$ occuring in action angle coordinates $(j_{i},\vartheta _{i})$. We show that, for a globally Hamiltonian vector
field $X_{f}$, 
\begin{equation}
\mathrm{\hat{e}}_{\ast }^{tZ_{f}}\sigma =\mathrm{e}^{-2\pi i\,tf/h}\,\mathrm{%
\widehat{e}}_{\ast }^{\,t \, \mathrm{lift}X_{f}}\sigma .  \label{1.2}
\end{equation}%
Replacing $f$ by a multivalued function $\vartheta $, defined up to an
arbitrary integer $n$, yields the multivalued expression 
\begin{equation}
\mathrm{\widehat{e}}_{\ast }^{\,tZ_{\vartheta }}\sigma =\mathrm{e}^{-2\pi
i\,t\vartheta /h}\,\mathrm{\widehat{e}}_{\ast }^{\,t\,\mathrm{lift}%
X_{\vartheta }}\sigma .  \label{1.3}
\end{equation}%
We observe that, for $t=h$, equation (\ref{1.3}) gives a single valued
expression 
\begin{equation}
\mathrm{\hat{e}}_{\ast }^{hZ_{\vartheta }}\sigma =\mathrm{e}^{-2\pi
i\,\vartheta }\,\mathrm{\widehat{e}}_{\ast }^{\,h\,\mathrm{lift}X_{\vartheta
}}\sigma (p).  \label{1.4}
\end{equation}%
The shifting operator 
\begin{equation}
\boldsymbol{a}_{X_{\vartheta }}=\mathrm{\widehat{e}}_{\ast
}^{\,hZ_{\vartheta }}=\mathrm{e}^{-2\pi i\,\vartheta }\,\mathrm{\widehat{e}}%
_{\ast }^{\,h\,\mathrm{lift}X_{\vartheta }}  \label{1.5}
\end{equation}%
is a skew adjoint operator on $\mathfrak{H,}$ which shifts the support of $%
\sigma \in \mathfrak{H}$ by $h$ in direction of $X_{\vartheta }$. If the
vector field $X_{\vartheta }$ is complete, then $\boldsymbol{a}%
_{X_{\vartheta }}^{n}=\mathrm{\widehat{e}}_{\ast }^{\,nhZ_{\vartheta }}$ for
every $n\in \mathbb{Z}$. If $\theta =2\pi \vartheta $ is a classical angle,
defined up to $2\pi n$, then $X_{\theta }=2\pi X_{\vartheta }$ and $h\,%
\mathrm{lift}X_{\vartheta }=\hbar \,\mathrm{lift}X_{\theta }$, where $\hbar
=h/2\pi .$ Therefore, we can write 
\begin{equation}
{\mathbf{a}}_{X_{\theta }}={\mathbf{a}}_{X_{2\pi \vartheta }}=
\mathrm{e}^{-i\theta }\,\mathrm{\widehat{e}}_{\ast }^{\,\hbar \,\mathrm{lift}X_{\theta}}.  \label{1.6}
\end{equation}

Our results provide an answer to Heisenberg's criticism that in 
Bohr-Sommerfeld theory there are not enough operators to describe
transitions between states \cite{heisenberg25}. In fact, these operators
exist, but it took quite a while to find them. \medskip 

In order to make the paper more accessible to the reader, we have provided
an introductory section with a comprehensive review of geometric quantization. 
Experts may omit this section and proceed directly to the next on Bohr-Sommerfeld theory. \medskip 

We would like to thank the referees for their comments and constructive
criticisms of an earlier version of this paper.

\section{Elements of geometric quantization}

%%%%%%%%%%%%%%%%

Let $(P, \omega)$ be a symplectic manifold. Geometric quantization can be
divided into three steps: prequantization, polarization, and unitarization.

\subsection{Principal line bundles with a connection}
%%%%%%%%%%%%%%%%%%

We begin with a brief review of connections on complex line bundles. \medskip

Let ${\mathbb{C} }^{\times }$ denote the mulitiplicative group of nonzero
complex numbers. Its Lie algebra ${\mathfrak{c}}^{\times }$ is isomorphic to
the abelian Lie algebra $\mathbb{C} $ of complex numbers. Different choices
of the isomorphism $\iota : \mathbb{C} \rightarrow {\mathfrak{c}}^{\times }$
lead to different factors in various expressions. Here to each $c \in 
\mathbb{C}$ we associate the $1$-parameter subgroup $t \mapsto {\mathrm{e}}%
^{2\pi i \, tc}$ of ${\mathbb{C} }^{\times }$. In other words, we take 
\begin{equation}
\iota : \mathbb{C} \rightarrow {\mathfrak{c}}^{\times }: c \mapsto \iota (c)
= 2\pi i \, c.  \label{eq-s2ss1newtwo}
\end{equation}

The prequantization structure for $(P,\omega )$ consists of a principal ${%
\mathbb{C} }^{\times }$ bundle ${\pi }^{\times }: L^{\times } \rightarrow P$
and a ${\mathfrak{c}}^{\times }$-valued ${\mathbb{C}}^{\times}$-invariant
connection $1$-form $\beta $ satisfying 
\begin{equation}
\mathrm{d} \beta = ({\pi }^{\times })^{\ast }(-%
\mbox{${\scriptstyle
\frac{{1}}{{h}}}$}\, \omega ),  \label{eq-s2ss1newthree}
\end{equation}
where $h$ is Planck's constant. The \emph{prequantization condition}
requires that the cohomology class $[-\frac{1}{h} \, \omega ]$ is integral,
that is, it lies in ${\mathrm{H}}^2(P, \mathbb{Z} )$. \medskip

Let $Y_c$ be the vector field on $L^{\times }$ generating the action of ${%
\mathrm{e}}^{2\pi i \, tc}$ on $L^{\times }$. In other words, the $1$%
-parameter group ${\mathrm{e}}^{t Y_c}$ of diffeomorphisms of $L^{\times }$
generated by $Y_c$ is 
\begin{equation}
{\mathrm{e}}^{t Y_c}: L^{\times } \rightarrow L^{\times }: {\ell }^{\times }
\rightarrow {\ell }^{\times } {\mathrm{e}}^{2\pi i \, tc} .
\label{eq-s2ss1newfour}
\end{equation}
The connection $1$-form $\beta $ is normalized by the requirement 
\begin{equation}
\langle \beta | Y_c \rangle =c.  \label{eq-s2ss1newfive}
\end{equation}
For each $c \ne 0$, the vector field $Y_c$ spans the vertical distribution $%
\mathrm{ver}\, TL^{\times }$ tangent to the fibers of ${\pi }^{\times }
:L^{\times } \rightarrow P$. The horizontal distribution $\mathrm{hor} \,
TL^{\times }$ on $L^{\times }$ is the kernel of the connection $1$-form $%
\beta $, that is, 
\begin{equation}
\mathrm{hor}\, TL^{\times } = \ker \beta .  \label{eq-s2ss1newsix}
\end{equation}
The vertical and horizontal distributions of $L^{\times }$ give rise to the
direct sum $TL^{\times } = \mathrm{ver}\, TL^{\times} \oplus \mathrm{hor}\,
TL^{\times }$, that is used to decompose any vector field $Z$
on $L^{\times }$ into its vertical and horizontal components, $Z = \mathrm{%
ver}\, Z + \mathrm{hor}\, Z$. Here the vertical component $\mathrm{ver}\, Z$
has range in $\mathrm{ver}\, TL^{\times }$ and the horizontal component has
range in $\mathrm{hor}\, TL^{\times }$. \medskip

If $X$ is a vector field on $P$, the unique horizontal vector field on $%
L^{\times }$, which is ${\pi }^{\times }$-related to $X$ is called the \emph{%
horizontal lift} of $X$ and is denoted by $\mathrm{lift}\, X$. In other
words, $\mathrm{lift}\, X$ has range in the horizontal distribution $\mathrm{%
hor}\, TL^{\times }$ and satisfies 
\begin{equation}
T{\pi }^{\times } \, \raisebox{2pt}{$\scriptstyle\circ \, $} \mathrm{lift}\,
X = X \, \raisebox{2pt}{$\scriptstyle\circ \, $} {\pi }^{\times }.
\label{eq-s2ss1eight}
\end{equation}

\noindent \textbf{Claim 2.1} \textit{A vector field $Z$ on $L^{\times }$ is
invariant under the action of ${\mathbb{C} }^{\times }$ on $L^{\times }$ if
and only if the horizontal component of $Z$ is the horizontal lift of its
projection $X$ to $P$, that is, $\mathrm{hor}\, Z = \mathrm{lift}\, X$ and
there is a smooth function $\kappa : P \rightarrow \mathbb{C} $ such that $%
\mathrm{ver}\, Z = Y_{\kappa (p)}$ on $L^{\times }_p = ({\pi }^{\times
})^{-1}(p)$.} \medskip

\noindent \textbf{Proof.} Since the direct sum $TL^{\times } = \mathrm{ver}%
\, TL^{\times} \oplus \mathrm{hor}\, TL^{\times }$ is invariant under the ${%
\mathbb{C} }^{\times}$ action on $L^{\times}$, it follows that the vector
field $Z$ is invariant under the action of ${\mathbb{C}}^{\times }$ if and
only if $\mathrm{hor}\, Z$ and $\mathrm{ver}\, Z$ are ${\mathbb{C} }^{\times
}$-invariant. But $\mathrm{hor}\, Z$ is ${\mathbb{C} }^{\times }$ invariant
if $T{\pi }^{\times } \, \raisebox{2pt}{$\scriptstyle\circ \, $} \mathrm{hor}%
\, Z = X \, \raisebox{2pt}{$\scriptstyle\circ \, $} {\pi }^{\times }$ for
some vector field $X$ on $P$, that is, $\mathrm{hor}\, Z = \mathrm{lift}\, X$%
. But this holds by definition. On the other hand, the vertical distribution 
$\mathrm{ver}\, TL^{\times }$ is spanned by the vector fields $Y_c$ for $c
\in \mathbb{C}$. Hence $\mathrm{ver}\, Z$ is ${\mathbb{C} }^{\times }$%
-invariant if and only if for every fiber $L^{\times}_p$ the restriction of $%
\mathrm{ver}\, Z$ to $L^{\times }_p$ coincides with the restriction of $Y_c$
to $L^{\times }_p$ for some $c \in \mathbb{C} $, that is, there is a smooth
complex valued function $\kappa $ on $P$ such that $c =\kappa (p)$. \hfill 
{\tiny $\blacksquare$} \medskip

Let $U$ be an open subset of $P$. A local smooth section $\tau :U\subseteq
P\rightarrow L^{\times }$ of the bundle ${\pi }^{\times }:L^{\times
}\rightarrow P$ gives rise to a diffeomorphism 
\begin{equation*}
{\eta }_{\tau }:L_{|U}^{\times }=\bigcup\limits_{p\in U}\big(({\pi }%
^{\times })^{-1}(p)\big)\rightarrow U\times {\mathbb{C}}^{\times }:{\ell }%
^{\times }\mapsto ({\pi }^{\times }({\ell }^{\times }),b)=(p,b),
\end{equation*}%
where $b\in {\mathbb{C}}^{\times }$ is the unique complex number such that ${%
\ell }^{\times }=\tau (p)b$. In the general theory of principal bundles the
structure group of the principal bundle acts on the right. In the theory of $%
{\mathbb{C}}^{\times }$ principal bundles, elements of $L^{\times }$ are
considered to be $1$-dimensional frames, which are usually written on the
right, see \cite{kostant}. The diffeomorphism ${\eta }_{\tau }$ is called a 
\emph{trivialization} of $L^{\times }_{\mid {U}}$. It intertwines the action of ${%
\mathbb{C}}^{\times }$ on the principal bundle $L^{\times }$ with the right
action of ${\mathbb{C}}^{\times }$ on $U\times {\mathbb{C}}^{\times }$,
given by mulitiplication in ${\mathbb{C}}^{\times }$. If a local section $%
\sigma :U\rightarrow L$ of $\pi :L\rightarrow P$ is nowhere zero, then it
determines a trivialization ${\eta }_{\tau }:L^{\times }_{\mid {U}}\rightarrow
U\times {\mathbb{C}}^{\times }$. Conversely, a local smooth section $\tau $
such that ${\eta }_{\tau }$ is a trivialization of $L^{\times }$ may be
considered as a local nowhere zero section of $L$. \medskip

In particular, for every $c \in \mathbb{C}$, which is identified with the
Lie algebra ${\mathfrak{c}}^{\times }$ of ${\mathbb{C} }^{\times }$,
equation (\ref{eq-s2ss1newthree}) gives ${\mathrm{e}}^{t\, Y_c} \, %
\raisebox{2pt}{$\scriptstyle\circ \, $} \tau = {\mathrm{e}}^{2\pi i \, tc}
\, \tau $. Differentiating with respect to $t$ and then setting $t=0$ gives 
\begin{equation}
Y_c \, \raisebox{2pt}{$\scriptstyle\circ \, $} \tau = 2\pi i \,c \, \tau .
\label{eq-s2ss1newtwelve}
\end{equation}

For every smooth complex valued function $\kappa : P \rightarrow \mathbb{C} $
consider the vertical vector field $Y_{\kappa }$ such that $Y_{\kappa }({%
\ell }^{\times }) = Y_{\kappa ({\pi }^{\times }({\ell }^{\times }))}$ for
every ${\ell }^{\times } \in L^{\times }$. The vector field $Y_{\kappa }$ is
complete and the $1$-parameter group of diffeomorphisms it generates is 
\begin{equation*}
{\mathrm{e}}^{t\, Y_{\kappa }}: L^{\times } \rightarrow L^{\times }: {\ell }%
^{\times } \mapsto {\ell }^{\times } {\mathrm{e}}^{2\pi i \, t\kappa ({\pi }%
^{\times }({\ell }^{\times }))}.
\end{equation*}
For every smooth section $\tau $ of the bundle ${\pi }^{\times }$ we have ${%
\mathrm{e}}^{t\, Y_{\kappa }} \, \raisebox{2pt}{$\scriptstyle\circ \, $}
\tau = {\mathrm{e}}^{2\pi i\, t\kappa } \, \tau $ so that 
\begin{equation}
Y_{\kappa } \, \raisebox{2pt}{$\scriptstyle\circ \, $} \tau = 2\pi i \,
\kappa \, \tau .  \label{eq-s2ss1newfourteen}
\end{equation}

Let $X$ be a vector field on $P$ and let $\mathrm{lift}\, X$ be its
horizontal lift to $L^{\times }$. The local $1$-parameter group ${\mathrm{e}}^{t\, \mathrm{lift}\, X}$ of local
diffeomorphisms of $L^{\times }$ generated by $\mathrm{lift}\, X$ commutes
with the action of ${\mathbb{C}}^{\times }$ on $L^{\times}$. For every ${%
\ell }^{\times }$, ${\mathrm{e}}^{t \, \mathrm{lift}\, X}({\ell }^{\times })$
is called \emph{parallel transport} of ${\ell }^{\times }$ along the
integral curve of $X$ starting at $p = {\pi }^{\times }({\ell }^{\times })$.
For every $p \in P$ the map ${\mathrm{e}}^{t\, \mathrm{lift}\, X}$ sends the
fiber $L^{\times }_p$ to the fiber $L_{{\mathrm{e}}^{tX}(p)}$. \medskip

There are several equivalent definitions of covariant derivative of a smooth
section of the bundle ${\pi }^{\times }$ in the direction of a vector field $%
X$ on $P$. We use the following one. The \emph{covariant derivative} of the
smooth section $\tau $ of the bundle ${\pi }^{\times }: L^{\times }
\rightarrow P$ in the direction $X$ is 
\begin{equation}
{\nabla }_X\tau = 
\mbox{${\displaystyle \frac{\dee }{\dee t}}
\rule[-10pt]{.5pt}{25pt} \raisebox{-10pt}{$\, {\scriptstyle t=0}$}$} ({%
\mathrm{e}}^{t \, \mathrm{lift}\, X})^{\ast} \tau .
\label{eq-s2ss1newfifteen}
\end{equation}

\noindent \textbf{Claim 2.2} \textit{The covariant derivative of a smooth
local section of the bundle ${\pi }^{\times }: L^{\times } \rightarrow P$ in
the direction $X$ is given by} 
\begin{equation}
{\nabla }_X\tau = 2\pi i \langle {\tau }^{\ast }\beta | X \rangle \, \tau .
\label{eq-s2ss1newsixteen}
\end{equation}

\noindent \textbf{Proof.} For every $p \in P$ we have 
\begin{align}
{\nabla }_X\tau (p) & = 
\mbox{${\displaystyle \frac{\dee }{\dee t}}
\rule[-10pt]{.5pt}{25pt} \raisebox{-10pt}{$\, {\scriptstyle t=0}$}$} 
({\mathrm{e}}^{t\, \mathrm{lift}\, X})^{\ast }\tau (p) = 
\mbox{${\displaystyle \frac{\dee }{\dee t}}
\rule[-10pt]{.5pt}{25pt} \raisebox{-10pt}{$\, {\scriptstyle t=0}$}$} 
({\mathrm{e}}^{-t \, \mathrm{lift}\, X} \, 
\raisebox{2pt}{$\scriptstyle\circ
\, $} \tau \, \raisebox{2pt}{$\scriptstyle\circ \, $} {\mathrm{e}}^{t\,
X})(p)  \notag \\
& = -\mathrm{lift}\, X (\tau (p)) + T\tau (X(p))  \notag \\
& = -\mathrm{lift}\, X (\tau (p)) + \mathrm{hor}\, (T\tau )X(p) + 
\mathrm{ver}\, (T\tau )X(p)  \notag \\
& = \mathrm{ver}\, (T\tau )X(p) .  \notag
\end{align}
The definition of the connection $1$-form $\beta $ and equation (\ref%
{eq-s2ss1newfourteen}) yield 
\begin{equation*}
\mathrm{ver}\, (T\tau (X(p)) = Y_{\langle \beta | T\tau \, %
\raisebox{2pt}{$\scriptstyle\circ \, $} X \rangle }(\tau (p)) = 2\pi i \,
\langle \beta | T\tau \, \raisebox{2pt}{$\scriptstyle\circ \, $} X \rangle
\tau (p).
\end{equation*}
Hence 
\begin{equation}
{\nabla }_X \tau = 2\pi i \, \langle \beta | T\tau \, \raisebox{2pt}{$%
\scriptstyle\circ \, $} X\rangle \, \tau ,  \label{eq-s2ss1newseventeen}
\end{equation}
which is equivalent to equation (\ref{eq-s2ss1newsixteen}). \hfill {\tiny $%
\blacksquare $}

\subsection{Associated line bundles}
%%%%%%%%%%%%%%%%

The complex line bundle $\pi : L \rightarrow P$ associated to the 
${\mathbb{C}}^{\times }$ principal bundle ${\pi }^{\times }: L^{\times } \rightarrow P$
is defined in terms of the action of ${\mathbb{C} }^{\times }$ on 
$(L^{\times } \times \mathbb{C})$ given by 
\begin{equation}
\Phi : {\mathbb{C}}^{\times } \times (L^{\times } \times \mathbb{C} )
\rightarrow L^{\times } \times \mathbb{C}: \big(b, ({\ell }^{\times }, c) %
\big) \mapsto ({\ell }^{\times }b , b^{-1}c) .  \label{eq-s2ss2neweighteen}
\end{equation}
Since the action $\Phi $ is free and proper, its orbit space $L = (L^{\times
} \times \mathbb{C})/{\mathbb{C}}^{\times }$ is a smooth manifold. A point $%
\ell \in L$ is the ${\mathbb{C} }^{\times }$ orbit $[({\ell }^{\times }, c)]$
through $({\ell }^{\times },c) \in (L^{\times } \times \mathbb{C})$, namely, 
\begin{equation}
\ell = [({\ell }^{\times }, c)] = \{ ({\ell }^{\times }b, b^{-1}c) \in
L^{\times } \times \mathbb{C} \, \rule[-4pt]{.5pt}{13pt}\, \, b \in {\mathbb{%
C} }^{\times } \} .  \label{eq-s2ss2newnineteen}
\end{equation}
The left action of ${\mathbb{C} }^{\times }$ on $\mathbb{C}$ gives rise to
the left action 
\begin{equation*}
\widehat{\Phi }: {\mathbb{C} }^{\times } \times L \rightarrow L : \big( a, [(%
{\ell }^{\times },c)] \big) \mapsto [({\ell }^{\times }, ac)],
\end{equation*}
which is well defined because $[({\ell }^{\times }, ac)] = [({\ell }^{\times
}b,b^{-1}(ac))] = [({\ell }^{\times }b, a(b^{-1}c))]$ for every ${\ell }%
^{\times } \in L^{\times}$, every $a$, $b \in {\mathbb{C} }^{\times}$ and
every $c \in \mathbb{C} $. The projection map ${\pi }^{\times }: L^{\times }
\rightarrow P$ induces the projection map 
\begin{equation*}
\pi : L \rightarrow L/{\mathbb{C} }^{\times } = P : \ell = [({\ell }^{\times
},c)] \mapsto \pi (\ell ) = \pi ([({\ell }^{\times }, c)]) = {\pi }^{\times
}({\ell }^{\times}) .
\end{equation*}

\noindent \textbf{Claim 2.3} \textit{A local smooth section $\sigma : U
\rightarrow L$ of the complex line bundle $\pi : L \rightarrow P$
corresponds to a unique mapping ${\sigma }^{\sharp}:L^{\times }_{\mid U}
\rightarrow \mathbb{C} $ such that for every $p \in U$ and every ${\ell }%
^{\times } \in L^{\times }_p$ 
\begin{equation}
\sigma (p) = [({\ell }^{\times }, {\sigma}^{\sharp }({\ell }^{\times }))],
\label{eq-s2ss2newtwentytwo}
\end{equation}
which is ${\mathbb{C} }^{\times }$-equivariant, that is, ${\sigma }^{\sharp}(%
{\ell }^{\times }b) = b^{-1}{\sigma }^{\sharp}({\ell }^{\times})$.} \medskip

\noindent \textbf{Proof.} Given $p \in U$ there exists $({\ell }^{\times },
c) \in L^{\times } \times \mathbb{C}$ such that $\sigma (p) = [({\ell }%
^{\times }, c)]$. Since the action of ${\mathbb{C} }^{\times }$ on $%
L^{\times }_p$ is free and transitive, it follows that the ${\mathbb{C} }%
^{\times }$ orbit $\{ ({\ell }^{\times }b, b^{-1}c) \in L^{\times }_p \times 
\mathbb{C} \, \rule[-4pt]{.5pt}{13pt}\, \, b \in {\mathbb{C} }^{\times } \} $
is the graph of a smooth function from $L^{\times }_p $ to $\mathbb{C} $,
which we denote by ${\sigma }^{\sharp}_p$. In particular, $c = {\sigma }%
^{\sharp}_p({\ell }^{\times })$ so that $\sigma (p) = [({\ell }^{\times},c)] = 
[({\ell }^{\times }, {\sigma }^{\sharp}_p({\ell }^{\times }))]$. As $%
p $ varies over $U$ we get a map 
\begin{equation*}
{\sigma }^{\sharp }: L^{\times }_{\mid U} \rightarrow \mathbb{C} : {\ell }^{\times
} \mapsto {\sigma }^{\sharp }({\ell }^{\times}) = {\sigma }^{\sharp }_{{\pi }%
^{\times}({\ell }^{\times })} ({\ell }^{\times }),
\end{equation*}
which satisfies equation (\ref{eq-s2ss2newtwentytwo}). For every $b \in {%
\mathbb{C} }^{\times }$ equations (\ref{eq-s2ss2newnineteen}) and (\ref%
{eq-s2ss2newtwentytwo}) imply that 
\begin{equation*}
\sigma (p) = [({\ell }^{\times }, {\sigma }^{\sharp} ({\ell }^{\times }))] =
[( {\ell }^{\times }b , b^{-1} {\sigma }^{\sharp}({\ell }^{\times }))] = [({%
\ell }^{\times }b, {\sigma }^{\sharp}({\ell }^{\times }b))].
\end{equation*}
Hence ${\sigma }^{\sharp }({\ell }^{\times }b) = b^{-1}{\sigma }^{\sharp }({%
\ell }^{\times })$. Thus the function ${\sigma }^{\sharp}$ is ${\mathbb{C} }%
^{\times }$-equivariant. \hfill {\tiny $\blacksquare$} \medskip

If $\tau : U \rightarrow L^{\times }$ is a local smooth section of the
bundle ${\pi }^{\times }: L^{\times } \rightarrow P$, then for every $p \in
P $ we have $\sigma (p) = [(\tau (p), {\sigma }^{\sharp}(\tau (p)) )]$ or $%
\sigma = [(\tau, {\sigma }^{\sharp} \, 
\raisebox{2pt}{$\scriptstyle\circ \,
$} \tau ) ]$ suppressing the argument $p$. The function $\psi = {\sigma }%
^{\sharp} \, \raisebox{2pt}{$\scriptstyle\circ \, $} \tau : U \rightarrow 
\mathbb{C}$ is the coordinate representation of the section $\tau $ in terms
of the trivialization ${\eta }_{\tau }: L^{\times }_{\mid U} \rightarrow U \times 
\mathbb{C}$. \medskip

\vspace{-.15in}Let $Z$ be a ${\mathbb{C} }^{\times }$-invariant vector field
on $L^{\times }$. Then $Z$ is ${\pi }^{\times }$-related to a vector field $%
X $ on $P$, that is, $T{\pi }^{\times } \, 
\raisebox{2pt}{$\scriptstyle\circ
\, $} Z = X \, \raisebox{2pt}{$\scriptstyle\circ \, $} {\pi }^{\times}$. We
denote by ${\mathrm{e}}^{t X}$ and ${\mathrm{e}}^{t Z}$ the local $1$%
-parameter groups of local diffeomorphisms of $P$ and $L^{\times }$
generated by $X$ and $Z$, respectively. Because the vector fields $X$ and $Z$
are ${\pi }^{\times }$-related, we obtain ${\pi }^{\times} \, %
\raisebox{2pt}{$\scriptstyle\circ \, $} {\mathrm{e}}^{t\, Z} = {\mathrm{e}}%
^{t\, X} \, \raisebox{2pt}{$\scriptstyle\circ \, $} {\pi }^{\times }$. In
other words, the flow ${\mathrm{e}}^{t\, Z}$ of $Z$ covers the flow ${%
\mathrm{e}}^{t \, X}$ of $X$. The local group ${\mathrm{e}}^{t\, Z}$ of
automorphisms of the principal bundle $L^{\times }$ act on the associated
line bundle $L$ by 
\begin{equation}
{\widehat{\mathrm{e}}}^{\, t\, Z}: L \rightarrow L : \ell = [({\ell }%
^{\times }, c)] \mapsto[ ({\mathrm{e}}^{t\, Z}( {\ell }^{\times }), c)] ,
\label{eq-s2ss2newtwentyfive}
\end{equation}
which holds for all $\ell =[({\ell }^{\times }, c)]$ for which ${\mathrm{e}}%
^{t\, Z}({\ell }^{\times })$ is defined. \medskip

\noindent \textbf{Lemma 2.4} \textit{The map ${\widehat{\mathrm{e}}}^{\, t
\, Z}$ is a local $1$-parameter group of local automorphisms of the line
bundle $L$, which covers the local $1$-parameter group ${\mathrm{e}}^{t\, X}$
of the vector field $X$ on $P$.} \medskip

\noindent \textbf{Proof.} We compute. For $\ell = [({\ell }^{\times }, c)]
\in L$ we have 
\begin{align}
{\widehat{\mathrm{e}}}^{\, (t+s)Z}(\ell ) & = {\widehat{\mathrm{e}}}^{\,
(t+s)Z} ([({\ell }^{\times }, c)]) = [({\mathrm{e}}^{(t+s)Z}({\ell }^{\times
}), c)] = [({\mathrm{e}}^{t\, Z}({\mathrm{e}}^{s\, Z}({\ell }^{\times })) ,c
)]  \notag \\
& = {\widehat{\mathrm{e}}}^{\, t\, Z}([ ({\mathrm{e}}^{s\, Z}({\ell }%
^{\times }, c)] = {\widehat{\mathrm{e}}}^{\, t\, Z} \, \raisebox{2pt}{$%
\scriptstyle\circ \, $} {\widehat{\mathrm{e}}}^{\, s\, Z} ([({\ell }^{\times
} ,c)]) = {\widehat{\mathrm{e}}}^{\, t\, Z} \, \raisebox{2pt}{$\scriptstyle%
\circ \, $} {\widehat{\mathrm{e}}}^{\, s\, Z}(\ell ). \notag
\end{align}
Hence ${\widehat{\mathrm{e}}}^{\, t\, Z}$ is a local $1$-parameter group of
local diffeomorphisms. Moreover, 
\begin{equation*}
\pi \, \raisebox{2pt}{$\scriptstyle\circ \, $} {\widehat{\mathrm{e}}}^{\, t
\, Z}(\ell ) = \pi ([ ({\mathrm{e}}^{t\, Z}({\ell }^{\times}) ,c)]) = {\pi}%
^{\times }({\mathrm{e}}^{t\, Z}({\ell }^{\times})) = {\mathrm{e}}^{t \, X}({%
\pi }^{\times }({\ell }^{\times }));
\end{equation*}
while 
\begin{equation*}
{\mathrm{e}}^{t\, X} \, \raisebox{2pt}{$\scriptstyle\circ \, $} \pi (\ell )
= {\mathrm{e}}^{t\, X}(\pi ([({\ell }^{\times },c)])) = {\mathrm{e}}^{t\, X}(%
{\pi }^{\times }({\ell }^{\times })).
\end{equation*}
This shows that ${\mathrm{e}}^{\, t\, Z}$ covers ${\mathrm{e}}^{t\, X}$.
Finally, for every $\ell = [({\ell }^{\times }, c)] \in L$ and every $b \in {%
\mathbb{C} }^{\times }$ 
\begin{equation*}
{\widehat{\Phi }}_b({\widehat{\mathrm{e}}}^{\, t \, Z} (\ell )) = {\widehat{%
\Phi }}_b([ ( {\mathrm{e}}^{t\, Z}({\ell }^{\times }), c)] ) = [({\Phi }_b({%
\mathrm{e}}^{t\, Z}({\ell }^{\times })) ,c)] = [ ( {\mathrm{e}}^{t\, Z}({%
\Phi }_b({\ell }^{\times})), c) ] ,
\end{equation*}
since $Z$ is a ${\mathbb{C} }^{\times }$-invariant vector field on $%
L^{\times }$. Therefore 
\begin{equation*}
{\widehat{\Phi }}_b({\widehat{\mathrm{e}}}^{\, t \, Z}(\ell )) = {\widehat{%
\mathrm{e}}}^{\, t \, Z}( [ ({\Phi }_b({\ell }^{\times }), c) ] ) = {%
\widehat{\mathrm{e}}}^{\, t \, Z} \, \raisebox{2pt}{$\scriptstyle\circ \, $} 
{\widehat{\Phi }}_b( [ ({\ell }^{\times }, c) ] ) = {\widehat{\mathrm{e}}}%
^{\, t \, Z} \, \raisebox{2pt}{$\scriptstyle\circ \, $} {\widehat{\Phi }}_b
(\ell ).
\end{equation*}
This shows that ${\widehat{\mathrm{e}}}^{\, t \, Z}$ is a local group of
automorphisms of the line bundle $\pi : L \rightarrow P$. \hfill {\tiny $%
\blacksquare $} \medskip

If $Z = \mathrm{hor}\, X$, then ${\mathrm{e}}^{t\, \mathrm{lift}\, X}({\ell }%
^{\times })$ is parallel transport of ${\ell }^{\times }$ along the integral
curve ${\mathrm{e}}^{t \, X}(p)$ of $X$ starting at $p= {\pi }^{\times }({%
\ell }^{\times })$. Similarly, if $\ell =[({\ell }^{\times }, c)] \in L$,
then 
\begin{equation}
{\widehat{\mathrm{e}}}^{\, t \, \mathrm{lift}\, X}(\ell ) = [( {\mathrm{e}}%
^{t\, \mathrm{lift}\, X}({\ell }^{\times }) ,c )]
\label{eq-s2ss2newtwentyeight}
\end{equation}
is parallel transport of $\ell \in L$ along the integral curve ${\mathrm{e}}%
^{t \, X}(p)$ of $X$ starting at $p$. The covariant derivative of a section $%
\sigma $ of the bundle $\pi : L \rightarrow P$ in the direction of the
vector field $X$ on $P$ is 
\begin{equation}
{\nabla }_X \sigma = 
\mbox{${\displaystyle \frac{\dee }{\dee t}}
\rule[-10pt]{.5pt}{25pt} \raisebox{-10pt}{$\, {\scriptstyle t=0}$}$} ({%
\widehat{\mathrm{e}}}^{\, t \, \mathrm{lift}\, X})^{\ast} \sigma = 
\mbox{${\displaystyle \frac{\dee }{\dee t}}
\rule[-10pt]{.5pt}{25pt} \raisebox{-10pt}{$\, {\scriptstyle t=0}$}$} ( {%
\widehat{\mathrm{e}}}^{\, -t \, \mathrm{lift}\, X} \, \raisebox{2pt}{$%
\scriptstyle\circ \, $} \sigma \, \raisebox{2pt}{$\scriptstyle\circ \, $} {%
\mathrm{e}}^{t\, X} ) .  \label{eq-s2ss2newtwentynine}
\end{equation}
Since ${\widehat{\mathrm{e}}}^{\, - t\, \mathrm{lift}\, X}$ maps ${\pi }%
^{-1}({\mathrm{e}}^{t\, X})$ onto ${\pi }^{-1}(p)$, equations (\ref%
{eq-s2ss2newtwentyeight}) and (\ref{eq-s2ss2newtwentynine}) are consistent
with the definitions in \cite{kobayashi-nomizu}. \medskip

\noindent \textbf{Theorem 2.5} \textit{Let $\sigma $ be a smooth section of
the complex line bundle $\pi : L \rightarrow P$ and let $X$ be a vector
field on $P$. For every ${\ell }^{\times } \in L^{\times}$} 
\begin{equation}
{\nabla }_X\sigma ({\pi }^{\times }({\ell }^{\times })) = [ ( {\ell }%
^{\times }, L_{\mathrm{lift}\, X}\big({\sigma }^{\sharp}({\ell }^{\times }) %
\big) ) ].  \label{eq-s2ss2newthirty}
\end{equation}
Here $L_X$ is the Lie derivative with respect to the vector field $X$. \medskip 

\noindent \textbf{Proof.} Let $p = {\pi }^{\times }({\ell }^{\times })$.
Equation (\ref{eq-s2ss2newtwentynine}) yields 
\begin{equation*}
{\nabla }_X \sigma (p) = 
\mbox{${\displaystyle \frac{\dee }{\dee t}}
\rule[-10pt]{.5pt}{25pt} \raisebox{-10pt}{$\, {\scriptstyle t=0}$}$} \big({%
\widehat{\mathrm{e}}}^{\, - t \, \mathrm{lift}\, X} \, \raisebox{2pt}{$%
\scriptstyle\circ \, $} \sigma \, \raisebox{2pt}{$\scriptstyle\circ \, $} {%
\mathrm{e}}^{t \, X}(\sigma (p)) \big) .
\end{equation*}
Recall that $\sigma (p) = [({\ell }^{\times }, {\sigma }^{\sharp}({\ell }%
^{\times }))] $. Hence 
\begin{equation*}
\sigma ({\mathrm{e}}^{t \, X}(p)) = [ ({\mathrm{e}}^{t\, \mathrm{lift}\, X}({%
\ell }^{\times }), {\sigma }^{\sharp } \big( {\mathrm{e}}^{t\, \mathrm{lift}%
\, X}({\ell }^{\times } )\big) ) ] .
\end{equation*}
By equation (\ref{eq-s2ss2newtwentyeight}) 
\begin{align}
{\widehat{\mathrm{e}}}^{\, - t \, \mathrm{lift}\, X} \big( \sigma ({\mathrm{e%
}}^{t \, X}(p)) \big) & = {\widehat{\mathrm{e}}}^{\, - t \, \mathrm{lift}\,
X} [ ({\mathrm{e}}^{t \, \mathrm{lift}\, X} ({\ell }^{\times }), {\sigma}%
^{\sharp }\big( {\mathrm{e}}^{t \, \mathrm{lift}\, X} ({\ell }^{\times }) %
\big) )]  \notag \\
&\hspace{-1in} = [ ({\mathrm{e}}^{\, - t \, \mathrm{lift}\, X} \big( {%
\mathrm{e}}^{t \, \mathrm{lift}\, X} \big) ({\ell }^{\times }), {\sigma}%
^{\sharp }\big( {\mathrm{e}}^{t \, \mathrm{lift}\, X} ({\ell }^{\times }) %
\big) )] = [({\ell }^{\times }, {\sigma}^{\sharp }\big( {\mathrm{e}}^{t \, 
\mathrm{lift}\, X} ({\ell }^{\times }) \big) )]  \notag
\end{align}
Therefore 
\begin{align}
{\nabla }_X \sigma (p) & = 
\mbox{${\displaystyle \frac{\dee }{\dee t}}
\rule[-10pt]{.5pt}{25pt} \raisebox{-10pt}{$\, {\scriptstyle t=0}$}$} \hspace{%
-3pt} {\widehat{\mathrm{e}}}^ {\, - t \, \mathrm{lift}\, X} \big( \sigma ({%
\mathrm{e}}^{t \, X}(p)) \big) = 
\mbox{${\displaystyle \frac{\dee }{\dee t}}
\rule[-10pt]{.5pt}{25pt} \raisebox{-10pt}{$\, {\scriptstyle t=0}$}$} \hspace{%
-3pt} [({\ell }^{\times}, {\sigma}^{\sharp } \big( {\mathrm{e}}^{t \, 
\mathrm{lift}\, X} ({\ell }^{\times }) \big) ) ]  \notag \\
&= [( {\ell }^{\times }, 
\mbox{${\displaystyle \frac{\dee }{\dee t}}
\rule[-10pt]{.5pt}{25pt} \raisebox{-10pt}{$\, {\scriptstyle t=0}$}$} {\sigma 
}^{\sharp}({\mathrm{e}}^{t\, \mathrm{lift}\, X}({\ell }^{\times })) )] = [ ({%
\ell }^{\times }, L_{\mathrm{lift}\, X} \big( {\sigma }^{\sharp}({\ell }%
^{\times }) \big) ) ].  \tag*{\tiny $\blacksquare $}
\end{align}

\subsection{Prequantization}

%%%%%%%%%%%%%%%%%%%%%%

Let $\pi : L \rightarrow P$ be the complex line bundle associated to the ${%
\mathbb{C} }^{\times }$ principal bundle ${\pi }^{\times } : L^{\times }
\rightarrow P$. The space $S^{\infty}(L)$ of smooth sections of $\pi : L
\rightarrow P$ is the representation space of prequantization. Since ${%
\mathbb{C}}^{\times } \subseteq \mathbb{C} $, we may identify $L^{\times }$
with the complement of the zero section in $L$. With this identification if $%
\sigma : U \rightarrow L$ is a local smooth section of $\pi : L \rightarrow
P $, which is nowhere vanishing, then it is a section of the bundle ${\pi }%
^{\times }_{\mid {L^{\times }_{\mid U}}}: L^{\times}_{\mid U} \rightarrow U$. \medskip

\noindent \textbf{Theorem 2.6} \textit{A ${\mathbb{C} }^{\times}$-invariant
vector field $Z$ on $L^{\times }$ preserves the connection $1$-form $\beta $
on $L^{\times }$ if and only if there is a function $f \in C^{\infty}(P)$
such that 
\begin{equation}
Z = \mathrm{lift}\, X_f - Y_{f/h},  \label{eq-s2ss3newthirtyone}
\end{equation}
where $h$ is Planck's constant.} \medskip

\noindent \textbf{Proof.} The vector field $Z$ on $L^{\times}$ preserves the
connection $1$-form $L_Z \beta =0$, which is equivalent to 
\begin{equation}
Z \mbox{$\, \rule{8pt}{.5pt}\rule{.5pt}{6pt}\, \, $} \mathrm{d} \beta = - 
\mathrm{d} (Z \mbox{$\, \rule{8pt}{.5pt}\rule{.5pt}{6pt}\, \, $} \beta ).
\label{eq-s2ss3newthirtytwo}
\end{equation}
Since $\mathrm{hor}\, Z \mbox{$\, \rule{8pt}{.5pt}\rule{.5pt}{6pt}\, \, $}
\beta =0$, it follows that $Z 
\mbox{$\, \rule{8pt}{.5pt}\rule{.5pt}{6pt}\,
\, $} \beta = \mathrm{ver}\, Z 
\mbox{$\, \rule{8pt}{.5pt}\rule{.5pt}{6pt}\,
\, $} \beta $. The ${\mathbb{C} }^{\times }$-invariance of $Z$ and $\beta $
imply the ${\mathbb{C}}^{\times }$-invariance of $\mathrm{ver}\, Z 
\mbox{$\,
\rule{8pt}{.5pt}\rule{.5pt}{6pt}\, \, $} \beta $. Hence $\mathrm{ver}\, Z %
\mbox{$\, \rule{8pt}{.5pt}\rule{.5pt}{6pt}\, \, $} \beta $ pushes forward to
a function ${\pi }_{\ast }(\mathrm{ver}\, Z 
\mbox{$\,
\rule{8pt}{.5pt}\rule{.5pt}{6pt}\, \, $} \beta ) \in C^{\infty}(P)$. Thus
the right hand side of equation (\ref{eq-s2ss3newthirtytwo}) reads 
\begin{equation}
-\mathrm{d} (Z \mbox{$\, \rule{8pt}{.5pt}\rule{.5pt}{6pt}\, \, $} \beta ) =
- ({\pi }^{\times})^{\ast }\big( \mathrm{d} ({\pi }^{\times }_{\ast }( 
\mathrm{ver}\, Z \mbox{$\, \rule{8pt}{.5pt}\rule{.5pt}{6pt}\, \, $} \beta )) %
\big).  \label{eq-s2ss3newthirtythree}
\end{equation}
By definition $Y_c \mbox{$\, \rule{8pt}{.5pt}\rule{.5pt}{6pt}\, \, $} \beta
= c$, for every $c \in \mathfrak{c}$. This implies 
\begin{equation*}
Y_c \mbox{$\, \rule{8pt}{.5pt}\rule{.5pt}{6pt}\, \, $} \mathrm{d} \beta =
L_{Y_c}\beta - \mathrm{d} (Y_c 
\mbox{$\, \rule{8pt}{.5pt}\rule{.5pt}{6pt}\,
\, $} \beta ) = 0.
\end{equation*}
Thus the left hand side of equation (\ref{eq-s2ss3newthirtytwo}) reads 
\begin{equation}
Z \mbox{$\, \rule{8pt}{.5pt}\rule{.5pt}{6pt}\, \, $} \mathrm{d} \beta = 
\mathrm{hor}\, Z \mbox{$\, \rule{8pt}{.5pt}\rule{.5pt}{6pt}\, \, $} \mathrm{d%
} \beta .  \label{eq-s2ss3newthirtyfour}
\end{equation}
The quantization condition (\ref{eq-s2ss1newthree}) together with (\ref%
{eq-s2ss3newthirtytwo}), (\ref{eq-s2ss3newthirtythree}) and (\ref%
{eq-s2ss3newthirtyfour}) allow us to rewrite equation (\ref%
{eq-s2ss3newthirtytwo}) in the form 
\begin{equation}
\mathrm{lift}\, X \mbox{$\, \rule{8pt}{.5pt}\rule{.5pt}{6pt}\, \, $} \big( ({%
\pi }^{\times })^{\ast } (-\mbox{${\scriptstyle \frac{{1}}{{h}}}$} \omega ) %
\big) = ({\pi }^{\times })^{\ast }\big( \mathrm{d} ({\pi }_{\ast }(\mathrm{%
ver}\, Z \mbox{$\, \rule{8pt}{.5pt}\rule{.5pt}{6pt}\, \, $} \beta )) \big) .
\label{eq-s2ss3newthirtyfive}
\end{equation}
Equation (\ref{eq-s2ss3newthirtyfive}) shows that $X$ is the Hamiltonian
vector field of the smooth function 
\begin{equation}
f = -h\, {\pi }_{\ast }(\mathrm{ver}\, Z 
\mbox{$\,
\rule{8pt}{.5pt}\rule{.5pt}{6pt}\, \, $} \beta ))
\label{eq-s2ss3newthirtysix}
\end{equation}
on $P$. We write $X = X_f$. This implies that 
\begin{equation}
\mathrm{hor}\, Z = \mathrm{lift}\, X_f.  \label{eq-s2ss3newthirtyseven}
\end{equation}

We still have to determine the vertical component $\mathrm{ver}\,Z$ of the
vector field $Z$. For each ${\ell }^{\times }\in L^{\times }$ there is a $%
c\in \mathfrak{c}$ such that $\mathrm{ver}\,Z=Y_{c}$. Since $Y_{c}$ is
tangent to the fibers of the ${\mathbb{C}}^{\times }$ principal bundle ${\pi 
}^{\times }:L^{\times }\rightarrow P$, the element $c$ of $\mathfrak{c}$
depends only on ${\pi }^{\times }({\ell }^{\times })=p\in P$. Therefore 
\begin{equation*}
-({\pi }_{\ast }^{\times }(\mathrm{ver}\,Z%
\mbox{$\,
\rule{8pt}{.5pt}\rule{.5pt}{6pt}\, \, $}\beta ))({\ell }^{\times })=-({\pi }%
_{\ast }^{\times }(Y_{c(p)}\mbox{$\, \rule{8pt}{.5pt}\rule{.5pt}{6pt}\, \, $}%
\beta ))({\ell }^{\times })=-c(p)=f(p)/h
\end{equation*}%
by equation (\ref{eq-s2ss3newthirtysix}). In other words, for every point ${%
\ell }^{\times }\in L^{\times }$ we have $\mathrm{ver}\,Z({\ell }^{\times
})=-Y_{f(p)/h}({\ell }^{\times })$, where $p={\pi }^{\times }({\ell }%
^{\times })$. Thus we have shown that 
\begin{equation*}
Z_{f}=Z=\mathrm{lift}\,X_{f}-Y_{f/h}.
\end{equation*}%
Reversing the steps in the above argument proves the converse. \hfill {\tiny 
$\blacksquare $} \medskip

To each $f\in C^{\infty }(P)$, we associate a prequantization operator 
\begin{equation}
{\mathcal{P}}_{f}:S^{\infty }(L)\rightarrow S^{\infty }(L):\sigma \mapsto {%
\mathcal{P}}_{f}\sigma =i\hbar \,%
\mbox{${\displaystyle \frac{\dee }{\dee t}}
\rule[-10pt]{.5pt}{25pt} \raisebox{-10pt}{$\, {\scriptstyle t=0}$}$}({%
\widehat{\mathrm{e}}}^{\,t\,Z_{f}})_{\ast }\sigma ,  
\label{eq-s2verynew}
\end{equation}%
where ${\widehat{\mathrm{e}}}^{\,t\,Z_{f}}$ is the action of 
${\mathrm{e}}^{t\,Z_{f}}:L^{\times }\rightarrow L^{\times }$ on $L$, see 
(\ref{eq-s2ss2newtwentyeight}). Note that the definition of covariant derivative
in equation (\ref{eq-s2ss2newtwentynine}) is defined in terms of the pull
back $({\widehat{\mathrm{e}}}^{\,t\,Z_{f}})^{\ast }\sigma $ of the section 
$\sigma $ by ${\widehat{\mathrm{e}}}^{\,t\,Z_{f}}$, while the prequantization
operator in (\ref{eq-s2verynew}) is defined using the push forward 
$({\widehat{\mathrm{e}}}^{\,t\,Z_{f}})_{\ast }\sigma $ of $\sigma $ by ${%
\widehat{\mathrm{e}}}^{\,t\,Z_{f}}$. \medskip

\noindent \textbf{Theorem 2.7} \textit{For every $f \in C^{\infty}(P)$ and
each $\sigma \in S^{\infty}(L)$} 
\begin{equation}
{\mathcal{P}}_f \sigma = (-i\hbar \, {\nabla }_{X_f} +f ) \sigma .
\label{eq-s2ss3newforty}
\end{equation}

\noindent \textbf{Proof.} Since the horizontal distribution on $L^{\times }$
is ${\mathbb{C}}^{\times }$-invariant and the vector field $Y_{c}$ generates
multiplication on each fiber of ${\pi }^{\times }$ by ${\mathrm{e}}^{2\pi
i\,c}$, it follows that ${\mathrm{e}}^{t\,\mathrm{lift}\,X_{f}}\,{\mathrm{e}}%
^{t\,Y_{f/h}}={\mathrm{e}}^{t\,Y_{f/h}}\,{\mathrm{e}}^{t\,\mathrm{lift}%
\,X_{f}}$. Since $f$ is conctant along integral curves of $X_{f}$, 
\begin{align}
{\mathrm{e}}^{t\,Z_{f}}& ={\mathrm{e}}^{t(\mathrm{lift}\,X_{f}-Y_{f/h})}  =
{\mathrm{e}}^{t\,\mathrm{lift}\,X_{f}}\,{\mathrm{e}}^{-t\,Y_{f/h}} \notag \\
&= {\mathrm{e}}^{t\,\mathrm{lift}\,X_{f}}{\mathrm{e}}^{-2\pi it\,f/h}={\mathrm{e}}^{-2\pi
i\,tf/h}{\mathrm{e}}^{t\,\mathrm{lift}\,X_{f}},  
\label{Z}
\end{align}%
and 
\begin{align}
{\mathcal{P}}_{f}\sigma & =i\hbar \,%
\mbox{${\displaystyle \frac{\dee }{\dee t}}
\rule[-10pt]{.5pt}{25pt} \raisebox{-10pt}{$\, {\scriptstyle t=0}$}$}({%
\widehat{\mathrm{e}}}^{\,t\,Z_{f}})_{\ast }\sigma =i\hbar \,%
\mbox{${\displaystyle \frac{\dee }{\dee t}}
\rule[-10pt]{.5pt}{25pt} \raisebox{-10pt}{$\, {\scriptstyle t=0}$}$}({%
\mathrm{e}}^{t\,\mathrm{lift}\,X_{f}}\,{\mathrm{e}}^{t\,Y_{f}/h})_{\ast
}\sigma  \notag \\
& =i\hbar \,%
\mbox{${\displaystyle \frac{\dee }{\dee t}}
\rule[-10pt]{.5pt}{25pt} \raisebox{-10pt}{$\, {\scriptstyle t=0}$}$}({%
\widehat{\mathrm{e}}}^{\,t\,\mathrm{lift}\,X_{f}})_{\ast }\sigma +i\hbar \,%
\mbox{${\displaystyle \frac{\dee }{\dee t}}
\rule[-10pt]{.5pt}{25pt} \raisebox{-10pt}{$\, {\scriptstyle t=0}$}$}({%
\widehat{\mathrm{e}}}^{\,t\,Y_{-f/h}})_{\ast }\sigma . 
\label{eq-s2ss3newfortyone}
\end{align}%
Since $({\widehat{\mathrm{e}}}^{\,t\,\mathrm{lift}\,X_{f}})_{\ast }\sigma =({%
\widehat{\mathrm{e}}}^{\,-t\,\mathrm{lift}\,X_{f}})^{\ast }\sigma $ equation
(\ref{eq-s2ss2newtwentynine}) gives 
\begin{equation}
i\hbar \,%
\mbox{${\displaystyle \frac{\dee }{\dee t}}
\rule[-10pt]{.5pt}{25pt} \raisebox{-10pt}{$\, {\scriptstyle t=0}$}$}({%
\widehat{\mathrm{e}}}^{\,t\,\mathrm{lift}\,X_{f}})_{\ast }\sigma =i\hbar \,%
\mbox{${\displaystyle \frac{\dee }{\dee t}}
\rule[-10pt]{.5pt}{25pt} \raisebox{-10pt}{$\, {\scriptstyle t=0}$}$}({%
\widehat{\mathrm{e}}}^{\,-t\,\mathrm{lift}\,X_{f}})^{\ast }\sigma =-i\hbar \,%
{\nabla }_{X_{f}}\sigma .  \label{eq-s2ss3newfortytwo}
\end{equation}%
Since $\pi \,\raisebox{2pt}{$\scriptstyle\circ \, $}{\widetilde{\mathrm{e}}}%
^{\,t\,Y_{-f/h}}=\pi \,\raisebox{2pt}{$\scriptstyle\circ \, $}{\mathrm{id}}%
_{P}$, where ${\mathrm{id}}_{P}$ is the identity map on $P$, it follows that 
\begin{equation*}
({\widehat{\mathrm{e}}}^{\,tY_{-f/h}})_{\ast }\sigma ={\widehat{\mathrm{e}}}%
^{\,t\,Y_{-f/h}}\,\raisebox{2pt}{$\scriptstyle\circ \, $}\sigma \,%
\raisebox{2pt}{$\scriptstyle\circ \, $}{\mathrm{id}}_{P}={\widehat{\mathrm{e}%
}}^{\,t\,Y_{-f/h}}\,\raisebox{2pt}{$\scriptstyle\circ \, $}\sigma .
\end{equation*}%
Let $\tau :U\subseteq P\rightarrow L^{\times }$ be a smooth local section of 
${\pi }^{\times }:L^{\times }\rightarrow P$, then $\sigma =[(\tau ,{\sigma }%
^{\sharp }\,\raisebox{2pt}{$\scriptstyle\circ \, $}\tau )]$. Thus for every $%
p\in P$ 
\begin{align}
{\widehat{\mathrm{e}}}^{\,-t\,Y_{f/h}}\,%
\raisebox{2pt}{$\scriptstyle\circ \,
$}\sigma (p)& ={\widehat{\mathrm{e}}}^{\,-t\,Y_{f/h}}[(\tau (p),{\sigma }%
^{\sharp }(\tau (p)))]=[({\mathrm{e}}^{-t\,Y_{f/h}}(\tau (p)),{\sigma }%
^{\sharp }(\tau (p)))]  \notag \\
& =[(\tau (p){\mathrm{e}}^{-2\pi i\,tf(p)/h},{\sigma }^{\sharp }(\tau
(p)))]=[(\tau (p),{\mathrm{e}}^{-2\pi i\,tf(p)/h}{\sigma }^{\sharp }(\tau
(p)))],  \notag
\end{align}%
since $[({\ell }^{\times }b,c)]=[({\ell }^{\times }b,b^{-1}(bc))]=[({\ell }%
^{\times },bc)]$ for every ${\ell }^{\times }\in L^{\times }$, $b\in {%
\mathbb{C}}^{\times }$ and $c\in \mathbb{C}$. It follows that 
\begin{equation}
{\widehat{\mathrm{e}}}^{\,-t\,Y_{f/h}}\,%
\raisebox{2pt}{$\scriptstyle\circ \,
$}\sigma (p)=[(\tau (p),{\mathrm{e}}^{-2\pi i\,tf(p)/h}{\sigma }^{\sharp
}(\tau (p)))]={\mathrm{e}}^{-2\pi i\,tf(p)/h}\sigma (p).  \label{Z2}
\end{equation}%
Therefore, 
\begin{align}
({\widehat{\mathrm{e}}}^{\,t\,Z_{f}})_{\ast }\sigma &
=({\widehat{\mathrm{e}}}^{\,t\,(\mathrm{lift}\,X_{f}-Y_{f/h})})_{\ast }\sigma  \notag \\
& =({\widehat{\mathrm{e}}}^{\,t\,\mathrm{lift}\,X_{f}}{\widehat{\mathrm{e}}}^
{\,-t\,Y_{f/h}})_{\ast}\sigma ={\mathrm{e}}^{-2\pi i\,tf(p)/h}
\big( {\widehat{\mathrm{e}}}^{\,t\, \mathrm{lift}\,X_{f}}\big) _{\ast }\sigma .  \label{Z3}
\end{align}

Since,
\begin{equation}
i\hbar \, \mbox{${\displaystyle \frac{\dee }{\dee t}}
\rule[-10pt]{.5pt}{25pt} \raisebox{-10pt}{$\, {\scriptstyle t=0}$}$}
{\widehat{\mathrm{e}}}_{\ast }^{\,t\,Y_{-f/h}}\sigma =i\hbar \,
\mbox{${\displaystyle \frac{\dee }{\dee t}}
\rule[-10pt]{.5pt}{25pt} \raisebox{-10pt}{$\, {\scriptstyle t=0}$}$}
({\mathrm{e}}^{-2\pi i\,tf/h}\sigma )=i\hbar (-2\pi i\,f/h)\sigma =f\,\sigma
\label{eq-s2ss3newfortyfour}
\end{equation}%
equations (\ref{eq-s2ss3newfortyone}), (\ref{eq-s2ss3newfortytwo}) and (\ref%
{eq-s2ss3newfortyfour}) imply equation (\ref{eq-s2ss3newforty}).
\hfill {\tiny $\blacksquare $} \medskip

A Hermitian scalar product $\langle \, \, | \, \, \rangle $ on the fibers of 
$L$ that is invariant \linebreak under parallel transport gives rise to a
Hermitian scalar product on the space $S^{\infty}(L)$ of smooth sections of $%
L$. Since the dimension of $(P, \omega )$ is $2k$, the scalar product of the
smooth sections ${\sigma }_1$ and ${\sigma }_2$ of $L$ is 
\begin{equation}
( {\sigma }_1 | {\sigma }_2 ) = \int_P \langle {\sigma }_1 | {\sigma }_2
\rangle \, {\omega }^k .
\label{H}
\end{equation}
The completion of the space $S^{\infty}_c(L)$ of smooth sections of $L$ with
compact support with respect to the norm $\| \sigma \| = \sqrt{(\sigma |
\sigma )}$ is the Hilbert space $\mathfrak{H}$ of the prequantization 
representation. \medskip

\noindent \textbf{Claim 2.8} \textit{The prequantization operator ${\mathcal{%
P}}_f$ is a symmetric operator on the Hilbert space $\mathfrak{H}$ of square
integrable sections of the line bundle $\pi : L \rightarrow P$ and satisfies
Dirac's quantization commutation relations} 
\begin{equation}
[{\mathcal{P}}_{f}, {\mathcal{P}}_{g}] = i\hbar \, {\mathcal{P}}_{\{ f, g \}
}.  \label{eq-s2ss3newfortyseven}
\end{equation}
for every $f$, $g \in C^{\infty}(P)$. Morover, the operator ${\mathcal{P}}_f$
is self adjoint if the vector field $X_f$ on $(P,\omega )$ is complete.
\medskip

\noindent \textbf{Proof.} We only verify that the commutation relations (\ref%
{eq-s2ss3newfortyseven}) hold. Let $f$, $g\in C^{\infty }(P)$ and let $%
\sigma \in S^{\infty }(L)$. We compute. 
\begin{align}
\lbrack {\nabla }_{X_{f}}-\mbox{${\scriptstyle \frac{{i}}{{\hbar}}}$}f,{%
\nabla }_{X_{g}}-\mbox{${\scriptstyle \frac{{i}}{{\hbar}}}$}g]\sigma & =[{%
\nabla }_{X_{f}},{\nabla }_{X_{g}}]\sigma +%
\mbox{${\scriptstyle
\frac{{i}}{{\hbar}}}$}\big({\nabla }_{X_{f}}(g\sigma )-g{\nabla }%
_{X_{f}}\sigma \big)  \notag \\
& \hspace{0.5in}-\mbox{${\scriptstyle \frac{{i}}{{\hbar}}}$}\big({\nabla }%
_{X_{g}}(f\sigma )-f{\nabla }_{X_{g}}\sigma \big)  \notag \\
& =\big(\lbrack {\nabla }_{X_{f}},{\nabla }_{X_{g}}]+%
\mbox{${\scriptstyle
\frac{{i}}{{\hbar }}}$}(L_{X_{f}}g-L_{X_{g}}f)\big)\sigma  \notag
\end{align}%
The quantization condition 
\begin{equation*}
\lbrack {\nabla }_{X_{f}},{\nabla }_{X_{g}}]-{\nabla }_{[X_{f},X_{g}]}=-%
\mbox{${\scriptstyle \frac{{i}}{{\hbar}}}$}\omega (X_{f},X_{g})
\end{equation*}%
yields 
\begin{equation*}
\lbrack {\nabla }_{X_{f}}-\mbox{${\scriptstyle \frac{{i}}{{\hbar}}}$}f,{%
\nabla }_{X_{g}}-\mbox{${\scriptstyle \frac{{i}}{{\hbar}}}$}g]={\nabla }%
_{[X_{f},X_{g}]}-\mbox{${\scriptstyle \frac{{i}}{{\hbar}}}$}\omega
(X_{f},X_{g})+\mbox{${\scriptstyle \frac{{i}}{{\hbar}}}$}%
(L_{X_{f}}g-L_{X_{g}}f)
\end{equation*}%
But $\{f,g\}=L_{X_{g}}f=-\omega (X_{f},X_{g})$. So $L_{X_{f}}g-L_{X_{g}}f=%
\{g,f\}-\{f,g\}=-2\{f,g\}$. Since $X_{g}%
\mbox{$\,
\rule{8pt}{.5pt}\rule{.5pt}{6pt}\, \, $}\omega =-\mathrm{d}g$, it follows
that 
\begin{equation*}
\lbrack X_{f},X_{g}]\mbox{$\, \rule{8pt}{.5pt}\rule{.5pt}{6pt}\, \, $}\omega
=L_{X_{f}}X_{g}\mbox{$\, \rule{8pt}{.5pt}\rule{.5pt}{6pt}\, \, $}\omega
=-L_{X_{f}}\mathrm{d}g=-\mathrm{d}L_{X_{f}}g=\mathrm{d}\{f,g\}.
\end{equation*}%
Consequently, $[X_{f},X_{g}]=-X_{\{f,g\}}$. So 
\begin{equation*}
\lbrack {\nabla }_{X_{f}}-\mbox{${\scriptstyle \frac{{i}}{{\hbar}}}$}f,{%
\nabla }_{X_{g}}-\mbox{${\scriptstyle \frac{{i}}{{\hbar}}}$}g]={\nabla }%
_{X_{\{f,g\}}}-\mbox{${\scriptstyle \frac{{i}}{{\hbar }}}$}\{f,g\}.\quad %
\mbox{\tiny $\blacksquare $}
\end{equation*}

\subsection{Polarization}

%%%%%%%%%%%%%%%%% 

Prequantization is only the first step of geometric quantization. The
prequantization operators do not satisfy Heisenberg's uncertainty relations.
In the case of Lie groups, the prequantization representation fails to be
\linebreak irreducible. These apparently unrelated shortcommings are
resolved by the next step of geometric quantization: the introduction of a
polarization. \medskip

A complex distribution $F \subseteq T^{\mathbb{C}} = \mathbb{C} \otimes TP$
on a symplectic manifold $(P, \omega )$ is \emph{Lagrangian} if for every $p
\in P$, the restriction of the symplectic form ${\omega }_p$ to the subspace 
$F_p \subseteq T^{\mathbb{C} }_p P$ vanishes identically, and ${\mathrm{rank}%
}_{\mathbb{C}}F = \mbox{$\frac{\scriptstyle 1}{\scriptstyle 2}\,$} \dim P$.
If $F$ is a complex distribution on $P$, let $\overline{F}$ be its complex
conjugate. Let 
\begin{equation*}
D = F \cap \overline{F} \cap T^{\mathbb{C} }P \, \, \, \mathrm{and} \, \, \,
E = (F + \overline{F}) \cap T^{\mathbb{C} } P .
\end{equation*}
A \emph{polarization} of $(P, \omega )$ is an involutive complex Lagrangian
distribution $F$ on $P$ such that $D$ and $E$ are involutive distributions
on $P$. Let $C^{\infty}(P)_F$ be the space of smooth complex valued
functions of $P$ that are constant along $F$, that is, 
\begin{equation}
C^{\infty}(P)_{F} = \{ f \in C^{\infty}(P) \otimes P \, 
\rule[-4pt]{.5pt}{13pt}\, \, \langle \mathrm{d} f | u \rangle = 0 \, \, 
\mbox{for every $u
\in F$} \} .  \label{eq-s2ss4newfortyeight}
\end{equation}
The polarization $F$ is \emph{strongly admissible} if the spaces $P/D$ and $%
P/E$ of integral manifolds of $D$ and $E$, respectively, are smooth
manifolds and the natural projection $P/D \rightarrow P/E$ is a submersion.
A strongly admissible polarization $F$ is locally spanned by Hamiltonian
vector fields of functions in $C^{\infty}(P)_F$. A polarization $F$ is \emph{%
positive} if $i\, \omega (u, \overline{u}) \ge 0$ for every $u \in F$. A
positive polarization $F$ is \emph{semi-definite} if $\omega (u, \overline{u}%
) =0 $ for $u \in F$ implies that $u \in D^{\mathbb{C} }$. \medskip

Let $F$ be a strongly admissible polarization on $(P,\omega )$. The space $%
S_{F}^{\infty }(L)$ of smooth sections of $L$ that are covariantly constant
along $F$ is the \emph{quantum space of states }corresponding to the
polarization $F$. \medskip 

The space $C_{F}^{\infty }(P)$ of smooth functions on $P$, whose Hamiltonian
vector field preserves the polarization $F$, is a Poisson subalgebra of $%
C^{\infty }(P)$. Quantization in terms of the polarization $F$ leads to \emph{%
quantization map} $\mathcal{Q}$, which is the restriction of the \emph{%
prequantization map} 
\begin{equation*}
\mathcal{P}:C^{\infty }(P)\times S^{\infty }(L)\rightarrow S^{\infty
}(L):(f,\sigma )\mapsto {\mathcal{P}}_{f}\sigma =(-i\hbar \,{\nabla }%
_{X_{f}}+f)\sigma
\end{equation*}%
to the domain $C_{F}^{\infty }(P)\times S_{F}^{\infty }(L)\subseteq
C^{\infty }(P)\times S^{\infty }(L)$ and the codomain $S_{F}^{\infty
}(L)\subseteq S^{\infty }(L)$. In other words, 
\begin{equation}
\mathcal{Q}:C_{F}^{\infty }(P)\times S^{\infty }(L)\rightarrow S_{F}^{\infty}(L):
(f,\sigma )\mapsto {\mathcal{Q}}_{f}\sigma =(-i\hbar \,{\nabla }%
_{X_{f}}+f)\sigma .  \label{eq-s2ss4newfortynine}
\end{equation}%
Quantization in terms of positive strongly admissible polarizations such
that $E\cap \overline{E}=\{0\}$ lead to unitary representations. For other
types of polarizations unitarity may require additional structure.

\section{Bohr-Sommerfeld theory}
%%%%%%%%%%%%%%%%%%% 

\subsection{Historical background}
%%%%%%%%%%%%%%%%%%%%%

Consider the cotangent bundle $T^{\ast }Q$ of a manifold $Q$. Let ${\pi }_Q:
T^{\ast }Q \rightarrow Q$ be the cotangent bundle projection map. The
Liouville $1$-form ${\alpha }_Q$ on $T^{\ast }Q$ is defined as follows. For
each $q \in Q$, $p\in T^{\ast }_qQ$ and $u_p \in T_p(T^{\ast }Q)$, 
\begin{equation}
\langle {\alpha }_Q | u_p \rangle = \langle p | T{\pi}_Q(u_p) \rangle.
\label{eq-s3ss1newfifty}
\end{equation}
The exterior derivative of ${\alpha }_Q$ is the canonical symplectic form $%
\mathrm{d} {\alpha }_Q$ on $T^{\ast }Q$. \medskip

\vspace{-.2in}Let $\dim Q = k$. A Hamiltonian system on $(T^{\ast }Q, 
\mathrm{d} {\alpha }_Q)$ with Hamiltonian $H_0$ is \emph{completely
integrable} if there exists a collection of $k-1$ functions $H_1, \ldots ,
H_{k-1} \in C^{\infty}(T^{\ast }Q)$, which are integrals of $X_{H_0}$, that
is, $\{ H_0 , H_i \} =0$ for $i = 1, \ldots , k-1$, such that $\{ H_i , H_j
\} = 0$ for $i$, $j = 1, \ldots , k-1$. Assume that the functions $H_0,
\ldots , H_{k-1}$ are independent on a dense open subset of $T^{\ast }Q$.
For each $p \in T^{\ast }Q$, let $M_p$ be the orbit of the family of
Hamiltonian vector fields $\{ X_{H_0}, \ldots , X_{H_{k-1}} \} $ passing
through $p$. This orbit is the largest connected immersed submanifold in 
$T^{\ast }Q$ with tangent space $T_{p^{\prime }}(M_p)$ equal to 
${\mathrm{span}}_{\mathbb{R} } \{ X_{H_0}(p^{\prime }), \ldots ,
X_{H_{k-1}}(p^{\prime }) \} $. The integral curve 
$t \mapsto {\mathrm{e}}^{t\, X_{H_0}}(p)$ of $X_{H_0}$ starting at $p$ is contained in $M_p$. Hence knowledge of the family $\{ M_p \, \rule[-4pt]{.5pt}{13pt}\, p \in 
T^{\ast }Q \} $ of orbits provides information on the evolution of the Hamiltonian
system with Hamiltonian $H_0$. \medskip

Bohr-Sommerfeld theory \cite{bohr}, \cite{sommerfeld} asserts that the
quantum states of the completely integrable system $(H_0, \ldots , H_{k-1},
T^{\ast }Q, \mathrm{d} {\alpha }_Q )$ are concentrated on the orbits $M \in
\{ M_p \, \rule[-4pt]{.5pt}{13pt}\, p \in T^{\ast }Q \} $, which satisfy the
\medskip

\noindent \textbf{Bohr-Sommerfeld condition}: For every closed loop $\gamma
: S^1 \rightarrow M\subseteq T^{\ast }Q$ there exists an integer $n$ such
that 
\begin{equation}
\oint {\gamma }^{\ast }({\alpha }_Q) = n\, h,  \label{eq-s3ss1newfiftyone}
\end{equation}
where $h$ is Planck's constant. \medskip

This theory applied to the bounded states of the relativistic hydrogen atom
yields results that agree exactly with the experimental data \cite%
{sommerfeld}. \linebreak Attempts to apply Bohr-Sommerfeld theory to the
helium atom, which is not completely integrable, failed to provide useful
results. In his 1925 paper \cite{heisenberg25} Heisenberg criticized
Bohr-Sommerfeld theory for not providing transition operators between
different states. At present, the Bohr-Sommerfeld \linebreak theory is
remembered by physicists only for its agreement with the 
quasi-classical limit of Schr\"{o}dinger theory. Quantum chemists have never stopped 
using it to describe the spectra of molecules.

\subsection{Geometric quantization in a toric polarization}

%%%%%%%%%%%%%%%

In order to interpret Bohr-Sommerfeld theory in terms of geometric
quantization, we consider a set $P\subseteq T^{\ast }Q$ consisting of points 
$p \in T^{\ast }Q$ where $X_{H_0}(p), \ldots , X_{H_{k-1}}(p)$ are linearly
independent and the orbit $M_p$ of the family $\{ X_{H_0}, \ldots ,
X_{H_{k-1}} \} $ of Hamiltonian vector fields on $(T^{\ast }Q, \mathrm{d} {%
\alpha}_{T^{\ast }Q})$ is diffeomorphic to the $k$ torus ${\mathbb{T}}^k = {%
\mathbb{R} }^k/{\mathbb{Z}}^k$. We assume that $P$ is a $2k$-dimensional
smooth manifold and that the set $B = \{ M_p \, \rule[-4pt]{.5pt}{13pt}\, p
\in P \} $ is a quotient manifold of $P$ with smooth projection map $\rho :
P \rightarrow B$. This implies that the symplectic form $\mathrm{d} {\alpha }%
_Q$ on $T^{\ast }Q$ restricts to a symplectic form on $P$, which we denote
by $\omega $. Let $D$ be the distribution on $P$ spanned by the Hamiltonian
vector fields $X_{H_0}, \ldots , X_{H_{k-1}}$. Since $\{ H_i, H_j \} =0 $
for $i$, $j =0, 1, \ldots , k-1$, it follows that $D$ is an involutive
Lagrangian distribution on $(P, \omega )$. Moreover, $F = D^{\mathbb{C} }$
is a strongly admissible \emph{polarization} of $(P, \omega )$. \medskip

Since the symplectic form $\mathrm{d}{\alpha }_{Q}$ on $T^{\ast }Q$ is
exact, the prequantization line bundle 
\begin{equation*}
{\pi }^{\times }:L_{T^{\ast }Q}^{\times }={\mathbb{C}}^{\times }\times
T^{\ast }Q\rightarrow T^{\ast }Q:\big( b, (q,p) \big) \mapsto (q,p)
\end{equation*}%
is trivial and has a connection $1$-form ${\beta }_{Q}=%
\mbox{${\scriptstyle \frac{{1}}{{2\pi i}}}$}\,\frac{\mathrm{d}b}{b}+%
\mbox{${\scriptstyle \frac{{1}}{{h}}}$}{\alpha }_{Q}$. Let $L^{\times }$ be
the restriction of $L_{T^{\ast }Q}^{\times }$ to $P$ and let $\alpha $ be
the $1$-form on $P$, which is the restriction of ${\alpha }_{Q}$ to $P$,
that is, $\alpha ={{\alpha }_{Q}}_{\mid {P}}$. Then $L^{\times }={\mathbb{C}}%
^{\times }\times P$ is a principal ${\mathbb{C}}^{\times }$ bundle over $P$
with projection map 
\begin{displaymath}
{\pi }^{\times }:L^{\times }={\mathbb{C}}^{\times }\times P\rightarrow
P:(b,p)\mapsto p
\end{displaymath}
and connection $1$-form $\beta =\mbox{${\scriptstyle \frac{{1}}{{2\pi i}}}$}%
\,\frac{\mathrm{d}b}{b}+\mbox{${\scriptstyle \frac{{1}}{{h}}}$}\,\alpha $.
The complex line bundle $\pi :L=\mathbb{C}\times P\rightarrow P:(c,p)\mapsto p$ 
associated to the principal bundle ${\pi }^{\times }$ is also trivial. Prequantization of this system is obtained by adapting the results of section 2. \medskip

Since integral manifolds of the polarization $D$ are $k$-tori, we have to
determine which of them admit nonzero covariantly constant sections of $L$.
\medskip

\noindent \textbf{Theorem 3.3} \textit{An integral manifold $M$ of the
distribution $D$ admits a section of the complex line bundle $L$, which is
nowhere zero when restricted to $M$, if and only if it satisfies the
Bohr-Sommerfeld condition} (\ref{eq-s3ss1newfiftyone}). \medskip

\noindent \textbf{Proof.} Supose that an integral manifold $M$ of $D$ admits
a nowhere zero section of $L_{\mid M}$. Since $\sigma $ is nowhere zero, it is a
section of $L^{\times }_{\mid M}$. Let $\gamma :S^{1}\rightarrow M$ be a loop in 
$M$. For each $t\in S^{1}$ let $\dot{\gamma}(t)\in T_{\gamma (t)}M$ be the
tangent vector to $\gamma $ at $t$. Since $\sigma $ is covariantly constant
along $M$, claim 2.2 applied to the section 
\begin{equation*}
\sigma :M\rightarrow L_{|M}^{\times }=\mathbb{C}\times M:p\mapsto (b(p),p)
\end{equation*}%
gives 
\begin{equation*}
{\nabla }_{X(p)}{\sigma }(p)=2\pi i\,\langle {\sigma }^{\ast }({\beta })(p)|X(p)\rangle \, 
\sigma (p)=0
\end{equation*}%
for every $p\in P$ and every $X(p)\in T_{p}M$. Taking $p=\gamma (t)$ and $%
X(p)=\dot{\gamma}(t)$ gives 
\begin{equation}
2\pi i\,\langle {\sigma }^{\ast }\beta (\gamma (t))|\dot{\gamma}(t)\rangle
\,\sigma (\gamma (t))=0.  \label{eq-s3ss2newfiftyfive}
\end{equation}%
Since $\beta =\frac{1}{2\pi i}\,\frac{\mathrm{d}b}{b}+\frac{1}{h}\alpha $
and $(\sigma \,\raisebox{2pt}{$\scriptstyle\circ \, $}\gamma )(t)=(b(\gamma
(t),\gamma (t)))$ we get 
\begin{align}
2\pi i\,\langle {\sigma }^{\ast }\beta (\gamma (t))|\dot{\gamma}(t)\rangle &
=2\pi i\,\langle \beta (\sigma (\gamma (t)))|\dot{\gamma}(t)\rangle  \notag
\\
& =\frac{1}{b(\gamma (t))}\frac{\mathrm{d}b(\gamma (t))}{\mathrm{d}t}+\frac{%
2\pi i}{h}\,\langle \alpha |\dot{\gamma}(t)\rangle  \notag \\
& =\frac{\mathrm{d}}{\mathrm{d}t}\ln b(\gamma (t))+\frac{2\pi i}{h}\,\langle
\alpha (\gamma (t))|\dot{\gamma}(t)\rangle .  \notag
\end{align}%
Hence equation (\ref{eq-s3ss2newfiftyfive}) is equivalent to 
\begin{equation*}
\frac{\mathrm{d}}{\mathrm{d}t}\ln b(\gamma (t))+\frac{2\pi i}{h}\langle
\alpha (\gamma (t))|\dot{\gamma}(t)\rangle =0,
\end{equation*}%
which integrated from $0$ to $2\pi $ gives 
\begin{equation*}
\ln b(\gamma (2\pi ))-\ln b(\gamma (0))=-\frac{2\pi i}{h}\,\int_{0}^{2\pi
}\langle \alpha (\gamma (t))|\dot{\gamma}(t)\rangle \,\mathrm{d}t=-\frac{%
2\pi i}{h}\oint {\gamma }^{\ast }\alpha .
\end{equation*}%
If $\gamma $ bounds a surface $\Sigma \subseteq M$, then Stokes' theorem
together with equation (\ref{eq-s3ss1newfiftyone}) and the quantization
condition (\ref{eq-s2ss1newthree}) yield 
\begin{equation*}
-\frac{2\pi i}{h}\,\oint {\gamma }^{\ast }\alpha =-\frac{2\pi i}{h}%
\,\int_{\Sigma }\mathrm{d}\alpha =-\frac{2\pi i}{h}\int_{\Sigma }\omega =0,
\end{equation*}%
because $M$ is a Lagrangian submanifold of $(P,\omega )$. Thus $\ln b(\gamma
(2\pi ))=\ln b(\gamma (0))$, which implies that the nowhere zero section $%
\sigma $ is parallel along $\gamma $. If $\gamma $ does not bound a surface
in $M$, but does satisfy the Bohr-Sommerfeld condition 
$\oint {\gamma }^{\ast }{\alpha }_{Q}=nh$ (\ref{eq-s3ss1newfiftyone}) with 
${\alpha }_{Q}$ replaced by its pull back $\alpha $ to $P$, then 
\begin{equation*}
\ln \Big(\frac{b(\gamma (2\pi ))}{b(\gamma (0))}\Big)=-\frac{2\pi i}{h}\oint 
{\gamma }^{\ast }\alpha =-\frac{2\pi i}{h}\,nh=-2\pi i\,n,
\end{equation*}%
so that 
\begin{equation*}
\frac{b(\gamma (2\pi ))}{b(\gamma (0))}={\mathrm{e}}^{-2\pi i\,n}=1.
\end{equation*}%
Hence $b(\gamma (2\pi ))=b(\gamma (0))$ and the nowhere zero section $\sigma 
$ is parallel along $\gamma $. \hfill {\tiny $\blacksquare $} \medskip

Note that that the manifolds $M$ that satisfy
Bohr-Sommerfeld conditions (\ref{eq-s3ss1newfiftyone}) are $k$%
-dimensional toric submanifolds of $P.$ We call them \emph{Bohr-Sommerfeld tori}. 
Let $\mathfrak{B}$ be the set of Bohr-Sommerfeld tori in $M$. Since
Bohr-Sommerfeld tori have dimension $k=\frac{1}{2}\dim P$, there is no
non-zero smooth section $\sigma _{0}:P\rightarrow L$ that is covariantly
constant along $D$. However, for each Bohr-Sommerfeld torus $M$, theorem 3.1
guarantees the existence of a non-zero smooth section 
$\sigma _{M}:M\rightarrow L_{\mid M}$, where $L_{\mid M}$ denotes the restriction of $L$ to $M$. \medskip

Let $\mathcal{S}=\{M\}$ be the set of Bohr-Sommerfeld tori in $P$. For each 
$M\in \mathcal{S}$, there exists a non-zero section $\sigma _{M}$ of $L$
restricted to $M$ determined up to a factor in ${\mathbb{C}}^{\times }$.
The direct sum 
\begin{equation}
\mathfrak{S}=\bigoplus\limits_{M\in \mathcal{S}}\{{\mathbb{C}}\sigma _{M}\}
\label{HD}
\end{equation}%
is the the space of quantum states of the Bohr-Sommerfeld theory. Thus, each
Bohr-Sommerfeld torus $M$ represents a 1-dimensional subspace $\{{\mathbb{C}}%
\sigma _{M}\}$ of quantum states. Moreover, $\{{\mathbb{C}}\sigma _{M}\}\cap
\{{\mathbb{C}}\sigma _{M^{\prime }}\}=\{0\}$ if $M\neq M^{\prime }$ because
Bohr-Sommerfeld tori are mutually disjoint. Hence, the collection 
$\{\sigma _{M} \}$ is a basis of $\mathfrak{S}.$\medskip 

For our toral polarization $F=D^{\mathbb{C}}$, the space of smooth functions
on $P$ that are constant along $F$, see equation (\ref{eq-s2ss4newfortyeight}), 
is $C_{F}^{\infty }(P)={\rho }^{\ast }(C^{\infty }(B))$, see lemma A.3.
For each $f\in C_{F}^{\infty }(P)$, the Hamiltonian vector field $X_{f}$ is
in $D$, that is, ${\nabla }_{X_{f}}\sigma _{M}=0$ for every basic state $%
\sigma _{M}\in \mathfrak{S}$. Hence the prequantization and quantization
operators act on the basic states $\sigma _{M}\in \mathfrak{S}$ by
multiplication by $f$, that is, 
\begin{equation}
{\mathcal{Q}}_{f}\sigma _{M}={\mathcal{P}}_{f}\sigma _{M}=
f\,\sigma _{M}=f_{\mid M} \, \sigma _{M}. 
\label{eq-s3ss2newfiftynine}
\end{equation}%
Note that $f_{\mid M}$ is a constant because $f\in C_{F}^{\infty }(P)$. For
a general quantum state $\sigma =\sum_{M\in \mathcal{S}}c_{M}\sigma _{M}\in 
\mathfrak{S},$ 
\begin{equation*}
{\mathcal{Q}}_{f}\sigma ={\mathcal{Q}}_{f}\sum_{M\in \mathcal{S}}c_{M}\sigma
_{M}=\sum_{M\in \mathcal{S}}c_{M}{\mathcal{Q}}_{f}\sigma _{M}=\sum_{M\in 
\mathcal{S}}c_{M}f_{\mid M}\,\sigma _{M}.
\end{equation*}

We see that, for every function $f\in C^{\infty }(P)$, each basic quantum
state $\sigma _{M}$ is an eigenstate of $\mathcal{Q}_{f}$ corresponding to
the eigenvalue $f_{\mid M}$. Since eigenstates corresponding to different
eigenvalues of the same symmetric operator are mutually orthogonal, it follows that
the basis $\{\sigma _{M}\}$ of $\mathfrak{S}$ is orthogonal. This is the
only information we have about scalar product in $\mathfrak{S}$. Our results
do not depend on other details about the scalar product in $\mathfrak{S}$.

\subsection{Shifting operators}
%%%%%%%%%%%%%%%%%%

\subsubsection{The simplest case $P=T^{\ast }\mathbb{T}^{k}$}
%%%%%%%%%%%%%%%%%%

We begin by assuming that $P=T^{\ast }\mathbb{T}^{k}$ with canonical
coordinates $(\boldsymbol{p,q})=(p_{1},...,p_{k},q_{1},...,q_{k})$ where,
for each $i=1,...,k$, $q_{i}$ is the canonical angular coordinate on the 
$i^{\text{th }}$torus and $p_{i}$ is the conjugate momentum. The symplectic
form is 
\begin{equation*}
\omega =\mathrm{d}\big( \sum_{i=1}^{k}p_{i}\mathrm{d}{q}_{i} \big)
=\sum_{i=1}^{k}\mathrm{d}p_{i}\wedge \mathrm{d}{q}_{i}.
\end{equation*}%
In this case, action-angle angle coordinates $(\boldsymbol{j},\boldsymbol{%
\vartheta })=(j_{1},\ldots ,j_{k},{\vartheta }_{1},\ldots ,{\vartheta }_{k})$
are obtained by rescaling the canonical coordinates so that, for every $%
i=1,...,k,$ we have $j_{i}=2\pi p_{i}$ and $\vartheta _{i}=q_{i}/2\pi $. Moreover,
the angle coordinate $\vartheta _{i}:T^{\ast }\mathbb{T}^{k}\rightarrow 
\mathbb{T}=\mathbb{R}/\mathbb{Z}$ is interpreted as a multi-valued real
function, the symplectic form 
\begin{equation}
{\omega }=\sum_{i=1}^{k}\mathrm{d}j_{i}\wedge \mathrm{d}{\vartheta }_{i},
\label{N1}
\end{equation}%
and the toric polarization of $(P,\omega )$ is given by 
$D=\mathrm{span}~\left\{ \frac{\partial }{\partial \vartheta _{1}},\ldots ,
\frac{\partial }{\partial \vartheta _{1}}\right\} .$\medskip

In terms of action-angle coordinates the Bohr-Sommerfeld tori in 
$T^{\ast }\mathbb{T}^{k}$ are given by equation 
\begin{equation}
\boldsymbol{j}=(j_{1},...,j_{k})=(n_{1}h,...,n_{k}h)=\boldsymbol{n}h,
\label{N1a}
\end{equation}%
where $\boldsymbol{n}=(n_{1},...,n_{k})\in \mathbb{Z}^{k}$. For each 
$\boldsymbol{n}\in \mathbb{Z}^{k}$, we denote by 
$\mathbb{T}_{\boldsymbol{n}}^{k}$ the corresponding Bohr-Sommerfeld torus in 
$\mathfrak{B}$. If $\beta 
=\frac{1}{2\pi i}\,\frac{\mathrm{d}b}{b}+\frac{1}{h}\mathrm{d}(
\sum_{i=1}^{k}j_{i}\, \mathrm{d}{\vartheta }_{i}) {\ }$is the connection
form in the principal line bundle 
$L^{\times } = {\C}^{\times } \times {\mathbb{T}}^k_{\mathbf{n}} \rightarrow 
{\mathbb{T}}^k_{\mathbf{n}}$, then sections 
\begin{equation}
\sigma _{\boldsymbol{n}}:\mathbb{T}_{\boldsymbol{n}}^{k}\rightarrow
L^{\times }:(\vartheta _{1},...,\vartheta _{k})\mapsto \mathrm{e}^{-2\pi
i(n_{1}\vartheta _{1}+...+n_{k}\vartheta _{k})},  \label{N1b}
\end{equation}
form a basis in the space $\mathfrak{S}$ of quantum states. 
\medskip 

For each $i=1,...,k,$ the vector field $\frac{\partial }{\partial j_{i}}$ is
transverse to $D$ and $-\frac{\partial }{\partial j_{i}} \lefthook 
\omega =- \mathrm{d}{\vartheta }_{i}$, so that $-\frac{\partial }{\partial j_i}$ 
is the Hamiltonian vector field of $\vartheta _{i}$. We write 
$X_{i} = -\frac{\partial }{\partial j_{i}}= X_{\vartheta _{i}}$. 
Equation (\ref{Z}) in section 2.1, for $f=\vartheta _{i},$ is multi-valued
because the phase factor is multi-valued, and %
\begin{equation}
\mathrm{e}^{tZ_{\vartheta _{i}}}=\mathrm{e}^{-2\pi it\vartheta _{i}/h}%
\mathrm{e}^{t\, \mathrm{lift}X_{i}}.  
\label{N2a}
\end{equation}

\begin{description}
\item[Claim 3.4] If $t=h$, then 
\begin{equation}
\mathrm{e}^{hZ_{X_{i}}}=\mathrm{e}^{-2\pi i\vartheta _{i}}\mathrm{e}^{h \, 
\mathrm{lift}X_{i}}.  \label{N3a}
\end{equation}%
is well defined.
\end{description}

\noindent \textbf{Proof}. For every $i=1,...,k$, consider an open interval $%
(a_{i},b_{i})$ in $\mathbb{R}$ such that $0 < b_{i}-a_{i}<1$. Let 
\begin{equation}
W=\vartheta _{1}^{-1}(a_{1},b_{1})\cap \vartheta _{2}^{-1}(a_{2},b_{2})\cap
...\cap \vartheta _{k}^{-1}(a_{k},b_{k}).  \label{N4a}
\end{equation}%
Since the action-angle coordinates $(j_1, \ldots , j_k, {\vartheta }_1, \ldots , 
{\vartheta}_k)$ are continuous, $W$ is an open subset of $%
P$. Let $\theta _{i}$ be a unique representative of ${\vartheta _i}_{\mid W}$ with values in 
$(a_{i},b_{i})$. With this notation,%
\begin{equation}
{\omega }_{\mid W}=\sum_{i=1}^{k}\mathrm{d}j_{i\mid W}\wedge \mathrm{d}{%
\theta }_{i}.  \label{N5a}
\end{equation}%
\medskip The restriction to $W$ of the vector field $X_{\vartheta _{i}}$ is
the genuinely Hamiltonian vector field of $\theta _{i}$, namely, 
\begin{equation}
X_{ {{\vartheta }_i}_{\mid W}} =X_{\theta _{i}}.  
\label{N6a}
\end{equation}%
The vector field 
\begin{equation}
Z_{\theta _{i}}=\mathrm{lift}\,X_{\theta _{i}}-Y_{\theta _{i}/h}  \label{N6b}
\end{equation}%
is well defined. Equation (\ref{Z}) yields ${\mathrm{e}}^{t\,Z_{\theta
_{i}}}={\mathrm{e}}^{-2\pi i\,t\theta _{i}/h}{\mathrm{e}}^{t\,\mathrm{lift}%
\,X_{\theta _{i}}}$. Hence  
\begin{equation}
{\mathrm{e}}^{h\,Z_{\theta _{i}}}={\mathrm{e}}^{-2\pi i\,\theta _{i}}{%
\mathrm{e}}^{h\,\mathrm{lift}\,X_{\theta _{i}}}.\medskip  \label{N7a}
\end{equation}

If we make another choice of intervals $(a_{i}^{\prime },b_{i}^{\prime })$
in $\mathbb{R}$ such that $0< b_{i}^{\prime }-a_{i}^{\prime }<1$ and let 
$W^{\prime}=\cap _{i=1}^{k}\vartheta _{i}^{-1}(a_{i}^{\prime },b_{i}^{\prime })$. Then 
$\theta _{i}^{\prime }$ with values in $(a_{i}^{\prime },b_{i}^{\prime })$
differs from $\theta _{i}$ by an integer, so that $\theta _{i}^{\prime}= 
\theta _{i}+n_{i}$, and in $W\cap W^{\prime }$, we have 
\begin{equation*}
{\mathrm{e}}^{-2\pi i\,\theta _{i}^{\prime }}=
{\mathrm{e}}^{-2\pi i\,(\theta _{i}+n_{i})}={\mathrm{e}}^{-2\pi i\,\theta _{i}}.
\end{equation*}
Moreover, $X_{\theta _{i}\mid W\cap W^{\prime }}=
X_{\theta _{i}^{\prime}\mid W\cap W^{\prime }}=
{X_i}_{\mid {W\cap W^{\prime }}}$, so that 
\begin{equation*}
({\mathrm{e}}^{h\,Z_{i}})_{\mid {L^{\times}_{\mid {W\cap W^{\prime }}}}} = 
({\mathrm{e}}^{h\,Z_{\theta _{i}}})_{\mid {L^{\times}_{\mid {W\cap W^{\prime }}}}}=
({\mathrm{e}}^{h\,Z_{\theta _{i}^{\prime }}})_{\mid L^{\times}_{\mid {W\cap W^{\prime }} }}.
\end{equation*}
Since we can cover $P$ by open contractible sets defined in equation (\ref%
{N4a}), we conclude that $\mathrm{e}^{hZ_{X_{i}}}$ is well defined by
equation (\ref{N3a}) and depends only on the vector field $X_{i}.$ \hfill 
{\tiny $\blacksquare $} \medskip

Consequently, there exists a connection preserving automorphism 
${\mathbf{A}}_{X_{i}}:L^{\times }\rightarrow L^{\times }$ such that, if $l^{\times }\in
L_{\mid {W}^{\times }}$, where $W\subseteq P$ is given by equation (\ref{N4a}), then 
\begin{equation}
\boldsymbol{A}_{X_{i}}(l^{\times })={\mathrm{e}}^{h\,Z_{i}}(l^{\times }).
\label{N8}
\end{equation}

\begin{description}
\item[Claim 3.5] The connection preserving automorphism 
$\boldsymbol{A} _{X_{i}}:L^{\times }\rightarrow L^{\times }$, defined by equation (\ref{N8})
depends only on the vector field $X_{i}$ and not the original choice of
the action-angle coordinates.
\end{description}

\noindent \textbf{Proof}. If $(j_{1}^{\prime },\ldots ,j_{k}^{\prime },{\vartheta }_{1}^{\prime },\ldots , {\vartheta }_{k}^{\prime })$ is another set of action-angle coordinates then 
\begin{equation}
j_{i}=\sum_{l=1}^{k}a_{il}\, j_{l}^{\prime }\text{ \ and \ }\vartheta
_{i}=\sum_{l=1}^{k}b_{il}\, \vartheta _{l}^{\prime },  \label{N9c}
\end{equation}%
where the matrices $A=(a_{il})$ and $B=(b_{il})$ lie in $\mathrm{Sl}(k, \Z)$ and 
$B=(A^{-1})^{T}$. In the new coordinates,  
\begin{equation*}
X_{\vartheta _{i}}=-\frac{\partial }{\partial j_{i}}=-\sum_{l=1}^{k}a_{il}%
\frac{\partial }{\partial j_{l}^{\prime }}=-\sum_{l=1}^{k}b_{il}X_{\vartheta _{l}^{\prime }}=X_{(b_{i1}\vartheta _{1}^{\prime }+...+b_{ik}\vartheta _{k}^{\prime })}.
\end{equation*}%
Clearly, 
\begin{equation}
{\mathrm{e}}^{\,t\,\mathrm{lift}\,X_{\vartheta _{i}}}=
{\mathrm{e}}^{t\, \mathrm{lift}\,X_{(b_{i1}\vartheta _{1}^{\prime }+ \cdots+
b_{ik}\vartheta _{k}^{\prime })}}.  \label{N10c}
\end{equation}%
In order to compare the phase factor entering equation (\ref{N2a}), we
consider an open contracible set $W\subseteq P$. \ As before, for each $%
i=1,...,k,$ choose a single-valued representative $\theta _{i}^{\prime }$ of 
$({\vartheta}^{\prime}_i)_{\mid W}$. Then 
\begin{equation}
\theta _{i}=\sum_{j=1}^{k}b_{ij}(\theta _{j}^{\prime}+l_{j})
=\sum_{j=1}^{k}b_{ij}\theta _{j}^{\prime}+\sum_{j=1}^{k}b_{ij}l_{j}
=\sum_{j=1}^{k}b_{ij}\theta _{j}^{\prime }+l,
\label{N11b}
\end{equation}%
where each $l_{j}$ is an integer and thus $l=\sum_{j=1}^{k}b_{ij}l_{j}$ is also
an integer. Hence, 
\begin{equation}
{\mathrm{e}}^{-2\pi i\,{\theta }_{i}}={\mathrm{e}}^{-2\pi i\,(b_{i1}\theta
_{1}^{\prime }+...+b_{ik}\theta _{k}^{\prime }+l)}={\mathrm{e}}^{-2\pi
i\,(b_{i1}\theta _{1}^{\prime }+...+b_{ik}\theta _{k}^{\prime })},
\label{N12a}
\end{equation}%
where $b_{i1},...,b_{ik}$ are integers. Since $l$ is constant, 
\begin{align}
X_{ {{\vartheta}_{i}}_{\mid W}} & =X_{\theta _{i}} =X_{(b_{i1}\theta _{1}^{\prime}+ 
\ldots +b_{ik}\theta _{k}^{\prime }+l)} \notag \\
& =X_{(b_{i1}\theta _{1}^{\prime}+ \ldots +b_{ik}\theta _{k}^{\prime })}=
{X_{(b_{i1}\vartheta _{1}^{\prime }+...+b_{ik}\vartheta _{k}^{\prime })}}_{\mid W}.  \label{N13a}
\end{align}%
Therefore, 
\begin{equation}
{\mathrm{e}}^{h\,Z_{\theta _{i}}}={\mathrm{e}}^{-2\pi i\,\theta _{i}}{%
\mathrm{e}}^{h\,\mathrm{lift}\,X_{\theta _{i}}}={\mathrm{e}}^{-2\pi
i\,(b_{i1}\theta _{1}^{\prime }+...+b_{ik}\theta _{k}^{\prime })}
{\mathrm{e}}^{h\,\mathrm{lift}\,X_{(b_{i1}\vartheta _{1}^{\prime }+ \ldots 
+b_{ik}\vartheta _{k}^{\prime })}},  \label{N14a}
\end{equation}%
which shows that the automorphism 
${\mathbf{A}}_{X_{\vartheta _{i}}}:L^{\times }\rightarrow L^{\times }$ depends 
on the vector field $X_{\vartheta _{i}}$ and \emph{not} on the action angle coordinates in which it is computed. \hfill {\tiny $\blacksquare $} 

\begin{description}
\item[Claim 3.6] For each $i=1,...,k$, the symplectomorpism 
$\mathrm{e}^{hX_{i}}:P\rightarrow P,$ where $h$ is Planck's constant, preserves the set 
$\mathcal{B}$ of Bohr-Sommerfeld tori in $P$.
\end{description}

\noindent \textbf{Proof.} Since $X_{i}$ is complete, $\mathrm{e}^{tX_{i}}:P\rightarrow P$ is a 1-parameter group of symplectomorphisms of $%
(P,\omega )$. Hence, $\mathrm{e}^{hX_{i}}:P\rightarrow P$ is well defined.
By equation (\ref{N1a}), ${j_i}_{\mid {\mathbb{T}}_{\boldsymbol{n}}^{k}}=n_{i}h$
for every Bohr-Sommerfeld torus $\mathbb{T}_{\boldsymbol{n}}^{k}$, where $%
\boldsymbol{n}=(n_{1},...n_{k})$. $\medskip $

Since $X_{i}=-\frac{\partial }{\partial j_{i}}$, 
\begin{align}
L_{X_{i}}(j_{i}\, \dee \vartheta _{i}) &= X_{i} \lefthook \dee j_{i}\wedge 
\dee \vartheta _{i}+\dee (X_{i} \lefthook j_{i} \, \dee \vartheta _{i})=-\dee \vartheta _{i},
\notag  \\
L_{X_{i}}(j_{l}\, \dee \vartheta _{l}) &=X_{i} \lefthook \dee j_{l}\wedge 
\dee \vartheta _{l}+\dee(X_{i} \lefthook j_{l} \, \dee \vartheta _{l})=0\, \, \, 
\mbox{for $l\neq i$.} \notag 
\end{align}%
This implies that, for every $l\neq i$, $( \mathrm{e}^{tX_{i}})^{\ast }
(j_{l}\, \dee \vartheta _{l}) =j_{l}\, \dee \vartheta _{l}$ and 
$( \mathrm{e}^{tX_{i}})^{\ast }(j_{i}\, \dee \vartheta _{i}) =(j_{i}-t)\dee \vartheta _{i}$. Therefore, if $\boldsymbol{j}=\boldsymbol{n}h$, then $(\mathrm{e}^{hX_{i}})^{\ast }j_{l}=j_{l}=n_{l}$, if $l\neq i$, and $( 
\mathrm{e}^{tX_{i}})^{\ast }j_{i}=(j_{i}-h)=(n_{i}-1)h$. This implies
that $\mathrm{e}^{hX_{\vartheta _{i}}}(\mathbb{T}_{\boldsymbol{n}}^{k})$ is
a Bohr-Sommerfeld torus. \hfill {\tiny $\blacksquare $} \medskip 

\vspace{-.2in}We denote by $\widehat{\mathbf{A}}_{X_{i}}:L\rightarrow L$ the action of 
${\mathbf{A}}_{X_{i}}:L^{\times }\rightarrow L^{\times }$. The automorphism 
$\widehat{\mathbf{A}}_{X_{i}}$ acts on sections of $L$ by pull-back and push-forward, namely, 
\begin{equation}
\begin{array}{rl}
(\boldsymbol{\widehat{A}}_{X_{i}})_{\ast }\sigma &=
({\widehat{\mathrm{e}}}^{\,h\,Z_{i}})_{\ast }\sigma = 
{\widehat{\mathrm{e}}}^{\, -\, h\, Z_{i}}\comp 
\sigma \comp {\mathrm{e}}^{\,h\,X_{i}},  \\
\rowspace (\boldsymbol{\hat{A}}_{X_{i}})^{\ast }\sigma &=
({\widehat{\mathrm{e}}}^{\,h\,Z_{i}})^{\ast }\sigma =
{\widehat{\mathrm{e}}}^{\,h\,Z_{i}}\comp \sigma \comp {\mathrm{e}}^{\,-h\,X_{i}}. 
\end{array} 
\label{15b}
\end{equation}%
Since $\boldsymbol{A}_{X_{i}\text{ }}:L^{\times }\rightarrow L^{\times }$ is
a connection preserving automorphism, it follows that if $\sigma $ satisfies
the Bohr-Sommerfeld conditions, then 
$(\boldsymbol{\widehat{A}}_{X_{i}})_{\ast}\sigma $ and 
$(\boldsymbol{\widehat{A}}_{X_{i}})^{\ast }\sigma $ also satisfy the 
Bohr-Sommerfeld conditions. In other words, 
$(\boldsymbol{\widehat{A}} _{X_{i}})_{\ast }$ and 
$(\boldsymbol{\widehat{A}}_{X_{i}})^{\ast }\ $preserve
the space $\mathfrak{S}$ of quantum states. The \emph{shifting operators} 
$\boldsymbol{a}_{X_{i}}$ and $\boldsymbol{b}_{X_{i}},$ corresponding to $X_{i},$ are the restrictions to $\mathfrak{S}$ of $(\boldsymbol{\widehat{A}}_{X_{i}})_{\ast }$ and 
$(\boldsymbol{\widehat{A}}_{X_{i}})^{\ast },$ respectively. 
For every $\boldsymbol{n}=(n_{1},...,n_{k})\in \mathbb{Z}^{k},$ equations (%
\ref{N1b}) and (\ref{N3a}) yield
\begin{equation}
\begin{array}{rl}
\boldsymbol{a}_{X_{i}} \sigma _{\boldsymbol{n}} &=
{\widehat{\mathrm{e}}}^{-\,h\,Z_{i}}\comp \sigma _{\boldsymbol{n}}\comp  {\widehat{\mathrm{e}}}^{\,h\,X_{i}}=\sigma _{\boldsymbol{n}_{i}^{-}}= 
\mathrm{e}^{-2\pi i(\sum_{j \ne i}n_{j}\vartheta _{j} + (n_i-1){\vartheta }_i)} 
{\sigma }_{\mathbf{n}}
\label{N16a} \\
\boldsymbol{b}_{X_{i}}\sigma _{\boldsymbol{n}} &={\widehat{\mathrm{e%
}}}^{\,h\,Z_{i}}\comp \sigma _{\boldsymbol{n}}\comp 
{\widehat{\mathrm{e}}}^{\,-h\,X_{i}}=\sigma _{\boldsymbol{n}_{i}^{+}}=
\mathrm{e}^{-2\pi i (\sum_{j \ne i}n_{j}\vartheta _{j} + (n_{i}+1){\vartheta}_i)} 
{\sigma }_{\mathbf{n}} .
\end{array}
\end{equation}%
For each $i=1,...,k$, $\boldsymbol{a}_{X_{i}}\hspace{-1pt}\raisebox{-1pt}{$\comp $} 
\boldsymbol{b}_{X_{i}} =%
\boldsymbol{b}_{X_{i}} \hspace{-1pt}\raisebox{-1pt}{$\comp $} \boldsymbol{a}_{X_{i}}=
\mathrm{id}_{\mathfrak{S}}$. In addition, the operators $\boldsymbol{a}_{X_{i}}$, $%
\boldsymbol{b}_{X_{j}\text{ }},$ for $i,j=1,...,k$, generate an abelian
group $\mathfrak{A}$ of linear transformations of $\mathfrak{S}$ into itself, which acts
transitively on the space of $1$-dimensional subspaces of $\mathfrak{S}$. \medskip 

Given a non-zero section $\sigma \in \mathfrak{S}$ supported on a
Bohr-Sommerfeld torus, the family of sections 
\begin{equation}
\{ ( \boldsymbol{a}_{X_{k}\text{ }}^{n_{k}} \cdots \boldsymbol{a}_{X_{1}%
\text{ }}^{n_{1}}\sigma ) \in \mathfrak{S}\setrule \,  n_{1},...n_{k}\in 
\mathbb{Z} \}   \label{N17a}
\end{equation}%
is a linear basis of $\mathfrak{S},$ invariant under the action \ of $%
\mathfrak{A}.$ Since $\mathfrak{A}$ is abelian, there exists a positive,
definite Hermitian scalar product $\left\langle \cdot \mid \cdot
\right\rangle $ on $\mathfrak{S}$, which is invariant under the action of $%
\mathfrak{A}$, and such that the basis in (\ref{N17a}) is orthonormal. It is
defined up to a constant positive factor. The completion of $\mathfrak{S}$
with respect to this scalar product yields a Hilbert space $\mathfrak{H}$ of
quantum states in the Bohr-Sommerfeld quantization of 
$T^{\ast }\mathbb{T}^{k}$. 

\subsubsection{General case}
%%%%%%%%%%%%%%%%%

Let $(P,\omega )$ be a symplectic manifold with toroidal polarization $D$,
and a covering by domains of action-angle coordinates. If $U$ and 
$U^{\prime}$ are the \linebreak 
domain of angle-action coordinates $(\boldsymbol{j},\boldsymbol{\vartheta })
=(j_{1},...,j_{k},\vartheta _{1},...,\vartheta _{k})$ and $(\boldsymbol{j}^{^{\prime }},\boldsymbol{\vartheta }^{\prime})= 
(j_{1}^{\prime },...,j_{k}^{\prime },\vartheta _{1}^{\prime},\ldots ,\vartheta _{k}^{\prime })$, respectively, and $U\cap U^{\prime }\neq
\varnothing ,$ then in $U\cap U^{\prime }$ we have 
\begin{equation}
j_{i}=\sum_{l=1}^{k}a_{il}\, j_{l}^{\prime }\text{ \ and \ }\vartheta _{i}
=\sum_{l=1}^{k}b_{il}\, \vartheta _{l}^{\prime },\text{\ }  \label{N18a}
\end{equation}%
where the matrices $A=(a_{il})$ and $B=(b_{il})$ lie in $\mathrm{Sl}(k, \Z)$ and 
$B=(A^{-1})^{T}$. \medskip 

Consider a complete locally Hamiltonian vector field $X$ on $(P,\omega )$
such that, for each angle-action coordinates 
$(\boldsymbol{j},\boldsymbol{\vartheta })$ with domain $U$, 
\begin{equation}
( X \lefthook \omega)_{\mid U}=-\dee (\boldsymbol{c}\cdot 
\boldsymbol{\vartheta })=- \dee (c_{1}\vartheta _{1}+\ldots +c_{k}\vartheta _{k}),  \label{N19a}
\end{equation}%
for some $\boldsymbol{c}=(c_{1},...,c_{k})\in \mathbb{Z}^{k}$. Equation (\ref{N18a}) shows that in $U\cap U^{\prime }$, we have 
\begin{displaymath}
c_{1}\vartheta _{1}+ \ldots +c_{k}\vartheta _{k}= 
c_{1}^{\prime }\vartheta _{1}^{\prime}+\ldots +c_{k}^{\prime }\vartheta _{k}^{\prime },
\end{displaymath} 
where $c_{i}^{\prime}=\sum_{j=1}^{k}c_{j}b_{ji}\in \mathbb{Z}$, for $i=1,...,k.$ As in the
preceding section, equation (\ref{Z}) with $f=\boldsymbol{c}\cdot 
\boldsymbol{\vartheta =}c_{1}\vartheta _{1}+...+c_{k}\vartheta _{k}$, which is
multi-valued, gives
\begin{equation}
\mathrm{e}^{tZ_{\boldsymbol{c}\cdot \boldsymbol{\vartheta }}}= 
\mathrm{e}^{-2\pi i\, t\, \boldsymbol{c}\cdot \boldsymbol{\vartheta }/h}
\mathrm{e}^{t\, \mathrm{lift}X},  \label{N20a}
\end{equation}%
which is multivalued, because the phase factor is multi-valued. 
As before, if we set $t=h$, we would get a single-valued expression 
$\mathrm{e}^{hZ_{\boldsymbol{c}\cdot 
\boldsymbol{\vartheta }}}=\mathrm{e}^{-2\pi i\boldsymbol{c}\cdot \boldsymbol{%
\vartheta }}\mathrm{e}^{h\, \mathrm{lift}X}$ because $c_{1},...,c_{k}\in 
\mathbb{Z}$. This would work along all integral curves 
$t\mapsto \mathrm{e}^{t\, X}(x)$ for $t\in \lbrack 0,1],$ which are contained in 
$U$.\medskip 

Consider now the case when, for $x_{0}\in U$, $\mathrm{e}^{hX}(x)\in
U^{\prime }$ and there exists $t_{1}\in (0,h)$ such that $x_{1}=\mathrm{e}%
^{t_{1}X}(x_{0})\in U\cap U^{\prime }$, where $U$ and $U^{\prime }$ are
domains of action-angle variables $(\boldsymbol{j},\boldsymbol{\vartheta })$
and $(\boldsymbol{j}^{\prime },\boldsymbol{\vartheta }^{\prime }),$
respectively. Moreover, assume that $\mathrm{e}^{tX}(x_{0})\in U$ for $%
t\in \lbrack 0,t_{1}]$ and $\mathrm{e}^{tX}(x_{1})\in U^{\prime }$ for $t\in
\lbrack 0,h-t_{1}].$ Using the multi-valued notation, for $l^{\times }\in L_{x_{0}}^{\times }$, 
we write  
\begin{align}
\boldsymbol{A}_{X}(l^{\times }) & =\mathrm{e}^{(h-t_{1})
Z_{\boldsymbol{c}^{\prime }\cdot \boldsymbol{\vartheta }^{\prime }}}( 
\mathrm{e}^{t_{1}Z_{\boldsymbol{c}\cdot \boldsymbol{\vartheta }}}(l^{\times })) \notag \\
& =\mathrm{e}^{-2\pi i(h-t_{1})\boldsymbol{c}^{\prime }\cdot 
\boldsymbol{\vartheta }^{\prime }/h}\mathrm{e}^{(h-t_{1})\mathrm{lift}X}( 
\mathrm{e}^{-2\pi it_{1}\boldsymbol{c}\cdot \boldsymbol{\vartheta }/h}
\mathrm{e}^{t_{1}\, \mathrm{lift}X}(l^{\times }))   \label{N21a} \\
&=( \mathrm{e}^{-2\pi i(h-t_{1})\boldsymbol{c}^{\prime }\cdot 
\boldsymbol{\vartheta }^{\prime }/h}\mathrm{e}^{-2\pi it_{1}\boldsymbol{c}%
\cdot \boldsymbol{\vartheta }/h}) \mathrm{e}^{(h-t_{1})\mathrm{lift}X}
( \mathrm{e}^{t_{1}\mathrm{lift}X}(l^{\times }))   \notag \\
&=\mathrm{e}^{-2\pi it_{1}( \boldsymbol{c}\cdot \boldsymbol{\vartheta
-c}^{\prime }\cdot \boldsymbol{\vartheta }^{\prime }) /h}
\mathrm{e}^{-2\pi i\boldsymbol{c}^{\prime }\cdot \boldsymbol{\vartheta }^{\prime }}%
\mathrm{e}^{h\mathrm{lift}X}(l^{\times }).  \notag
\end{align}
Let $W$ be a neighbourhood of $x_{1}$ in $P$ such that $U\cap W$ and $%
U^{\prime }\cap W^{\prime }$ are contractible. For each $i=1,...,k,$ let $%
\theta _{i}$ be a single-valued representative of $\vartheta _{i}$ as in the
proof of claim 3.3.2. Similarly, we denote by $\theta _{i}^{\prime }$ a
single-valued representative of $\vartheta _{i}^{\prime }$. Equation 
(\ref{N19a}) shows that in $U\cap U^{\prime }\cap W$, the functions 
$c_{1}\theta _{1}+\cdots +c_{k}\theta _{k}$ and 
$c_{1}^{\prime }\theta _{1}^{\prime}+\cdots +c_{k}^{\prime }\theta _{k}^{\prime }$ are local Hamiltonians of the vector field $X$ and are constant along the integral curve of 
$X_{\mid W}$. Hence, we have to make the choice of repressentatives $\theta _{i}$ and 
$\theta _{i}^{\prime }$ so that 
\begin{equation}
c_{1}\theta _{1}(x_{1})+ \cdots +c_{k}\theta _{k}(x_{1})=
c_{1}\theta _{1}^{\prime}(x_{1})+\cdots +c_{k}\theta _{k}^{\prime }(x_{1}).  
\label{N22a}
\end{equation}%
With this choice, $\mathrm{e}^{-2\pi it_{1}( \boldsymbol{c}\cdot 
\boldsymbol{\vartheta -c}^{\prime }\cdot \boldsymbol{\vartheta }^{\prime }) /h}=1$, and 
\begin{equation}
\boldsymbol{A}_{X}(l^{\times })=\mathrm{e}^{-2\pi i\boldsymbol{c}^{\prime}\cdot \boldsymbol{\vartheta }^{\prime }}\mathrm{e}^{h\mathrm{lift}X}(l^{\times })  \label{N23a}
\end{equation}%
is well defined. It does not depend on the choice of the intermediate point $x_{1}$ in $U\cap U^{\prime }$.$\medskip $

In the case when $m+1$ action-angle coordinate charts with domains $%
U_{0},U_{1,}...,U_{m}$ are needed reach $x_{m}=%
\mathrm{e}^{hX}(x_{0})\in U_{m}$ from $x_{0}\in U_{0}$, we choose $x_{1}=
\mathrm{e}^{t_{1}X}(x_{0})\in U_{0}\cap U_{1},$ $x_{2}=\mathrm{e}^{t_{2}X}(x_{1})\in
U_{1}\cap U_{2}, \ldots , x_{m-1}=\mathrm{e}^{t_{m-1}X}(x_{m-2})\in U_{m-1}$
and end with $x_{m}=\mathrm{e}^{(h-t_{1}-\ldots -t_{m-1})X}(x_{m-1})\in U_{m}$.
At each intermediate point $x_{1}, \ldots ,x_{m-1},$ we repeat the the argument
of the preceding paragraph. We conclude that there is a well defined
connection preserving automorphism 
$\boldsymbol{A}_{X}:L^{\times}\rightarrow L^{\times }$ 
is well defined by the procedure given here, and
it depends only on the complete locally Hamiltonian vector field $X$
satisfying condition (\ref{N19a}). The automorphism 
$\boldsymbol{A} _{X}:L^{\times }\rightarrow L^{\times }$ of the principal bundle 
$L^{\times}$ leads to an automorphism 
$\widehat{A} _{X}$ of the associated line bundle $L$. As in equation 
(\ref{15b}), the shifting operators corresponded to the complete locally 
Hamiltonian vector field $X$ are 
\begin{equation}
\begin{array}{rl}
\boldsymbol{a}_{X} :\mathfrak{S}\rightarrow \mathfrak{S}:\sigma \mapsto 
(\boldsymbol{\widehat{A}}_{X})_{\ast }\sigma ,  \\
\rowspace \boldsymbol{b}_{X} : \mathfrak{S}\rightarrow \mathfrak{S}:\sigma \mapsto 
(\boldsymbol{\widehat{A}}_{X})^{\ast }\sigma .  
\end{array} 
\label{N24a}
\end{equation}

In absence of monodromy, if we have $k$ independent, complete, locally
Hamiltonian vector fields $X_{i}$ on $(P,\omega )$ that satisfy the
conditions leading to equation (\ref{N19a}), then the operators 
$\boldsymbol{a}_{X_{i}}$, $\boldsymbol{b}_{X_{j}}$ for $i,j=1,...,k$ generate an
abelian group $\mathfrak{A}$ of linear transformations of $\mathfrak{S}$. If
the local lattice $\mathfrak{B}$ of Bohr-Sommerfeld tori is regular, then 
$\mathfrak{A}$ acts transitively on the space of $1$-dimensional subspaces of 
$\mathfrak{S}$. This enables us to construct an $\mathfrak{A}$-invariant
Hermitian scalar product on $\mathfrak{S}$, which is unique up to an
arbitrary positive constant. The completion of $\mathfrak{S}$ with respect to
this scalar product yields a Hilbert space $\mathfrak{H}$ of quantum states
in the Bohr-Sommerfeld quantization of $(P,\omega )$.$\medskip $

\subsubsection{When things go wrong}
%%%%%%%%%%%%%%%%%%

\paragraph{Monodromy.}
%%%%%%%%%%%%%%%%%

In presence of monodromy, there may exist loops in the local lattice $\mathfrak{B}$ of
Bohr-Sommerfeld tori such that for some ${\alpha }_1, \ldots , {\alpha }_m 
\in \{ 1, \ldots , k \}$ the mapping 
\begin{equation*}
\big( e^{hX_{\alpha _{m}}}\comp  \cdots \comp 
{\mathrm{e}}^{hX_{\alpha _{1}}}\big)_{\mid {\mathbb{T}}_{\boldsymbol{n}}^{k}}:
\mathbb{T}_{\boldsymbol{n}}^{k}\rightarrow \mathbb{T}_{\boldsymbol{n}}^{k}
\end{equation*}%
is not the identity on $\mathbb{T}_{\boldsymbol{n}}^{k}$. In this case
shifting operators are multivalued, and there exists a phase factor 
$\mathrm{e}^{i\varphi }$ such that 
\begin{equation*}
( {\mathbf{a}}_{X_{\alpha _{m}}}\comp \cdots \comp 
{\mathbf{a}}_{X_{\alpha _{1}}}) \sigma _{\boldsymbol{n}}= 
\mathrm{e}^{i\varphi }\sigma _{\boldsymbol{n}}.
\end{equation*}%
Nevertheless, we still may use shifting operators to define a Hilbert space
structure in $\mathfrak{S}$. 

\paragraph{Incompleteness of $X$.}
%%%%%%%%%%%%%%%%%

If a locally Hamiltonian vector field $X$ on $(P,\omega ),$ which satisfies
the conditions leading to equation (\ref{N19a}),  is incomplete, then $%
\mathrm{e}^{hX}$ is not globally defined. If the integral curve $\mathrm{e}%
^{tX}(p)$ of $X$ originating at $p$ is defined for $t\in (t_{\min },t_{\max
}),$ then $\mathrm{e}^{hX}\left( \mathrm{e}^{tX}(p)\right) $ is defined for $%
t\in (t_{\min },t_{\max }-h),$ and $\mathrm{e}^{-hX}\left( \mathrm{e}%
^{tX}(p)\right) $ is defined for $t\in (t_{\min }+h,t_{\max })$.
Accordingly, the corresponding shifting operators $\boldsymbol{a}_{X}$ and $%
\boldsymbol{b}_{X}$ are not globally defined on the space of quantum states $%
\mathfrak{S}$. This usually occurs in systems, which do not have regular
toral polarization, and we consider only an open dense part of the phase
space where the toral polarization is regular. 

\subsection{Local lattice structure}

%%%%%%%%%%%%%%%

The discussion in section 3.2 did not address the question of labeling the
sections ${\sigma }_b$ in $\mathfrak{B}$ of the toral polarization $D$ by
the quantum numbers $\mathbf{n} = (n_1, \ldots , n_k)$ associated to the
Bohr-Sommerfeld $k$-torus $T = M_b$, which is the support of ${\sigma }_b$.
\medskip

These quantum numbers \emph{do} depend on the choice of action angle
coordinates. If $(j^{\prime }, {\vartheta }^{\prime }) \in V \times {\mathbb{%
T}}^k$ is another choice of action angle coordinates in the trivializing
chart $(U^{\prime }, {\psi }^{\prime })$, where $T \subseteq U^{\prime }$,
then the quantum numbers ${\mathbf{n}}^{\prime }$ of $T$ in $(j^{\prime }, {%
\vartheta }^{\prime })$ coordinates are related to the quantum numbers $%
\mathbf{n}$ of $T$ in $(j, \vartheta )$ coordinates by a matrix $A \in 
\mathrm{Gl} (k, \mathbb{Z} )$ such that ${\mathbf{n}}^{\prime }= A\, {%
\mathbf{n}}$, because by claim A.8 in the appendix on $U \cap U^{\prime }$
the action coordinates $j^{\prime }$ is related to the action coordinate $j$
by a constant matrix $A \in \mathrm{Gl} (k, \mathbb{Z})$. Let ${\mathbb{L}}%
_{\mid U} = \{ \mathbf{n} \in {\mathbb{Z} }^k \, \rule[-4pt]{.5pt}{13pt}\, T_{%
\mathbf{n}} \subseteq U\} $. Then ${\mathbb{L}}_{\mid U}$ is the \emph{local
lattice structure} of the Bohr-Sommerfeld tori $T_{\mathbf{n}}$, which lie
in the action angle chart $(U, \psi )$. If $(U, \psi )$ and $(U^{\prime }, {%
\psi }^{\prime })$ are action angle charts, then the set of Bohr-Sommerfeld
tori in $U\cap U^{\prime }$ are \emph{compatible}. More precisely, on $U
\cap U^{\prime }$ the local lattices ${\mathbb{L}}_{\mid U}$ and ${\mathbb{L}}%
_{\mid {U^{\prime }}}$ are compatible if there is a matrix $A \in \mathrm{Gl}(k, 
\mathbb{Z})$ such that ${\mathbb{L}}_{\mid {U^{\prime }}} = A \, {\mathbb{L}}_{\mid U}$.
Let $\mathcal{U} ={\ \{ U_i \} }_{i \in I}$ be a \emph{good} covering of $P$%
, that is, every finite intersection of elements of $\mathcal{U}$ is either
contractible or empty, such that for each $i \in I$ we have a trivializing
chart $(U_i , {\psi }_i)$ for action angle coordinates for the toral bundle 
$\rho : P \rightarrow B$. Then ${\{ {\mathbb{L}}_{U_u} \} }_{i \in I}$ is a
collection of pairwise compatible local lattice structures for the
collection $\mathcal{S}$ of Bohr-Sommerfeld tori on $P$. We say that 
$\mathcal{S}$ has a \emph{local lattice structure}. \medskip

The next result shows how the operator 
$({\widehat{\mathrm{e}}}^{\, h\, Z_{{\vartheta }_i}})_{\ast }$ of section 3.3 affects the quantum numbers of the Bohr-Sommerfeld torus $T = T_{\mathbf{n}}$. \medskip

\noindent \textbf{Claim 3.8} \textit{Let $(U , \psi )$ be a chart in $(P, \omega )$ for action angle coordinates $(j, \vartheta )$. For every
Bohr-Sommerfeld torus $T = T_{\mathbf{n}}$ in $U$ with quantum numbers 
$\mathbf{n}= (n_1, \ldots , n_k)$, the torus 
${\mathrm{e}}^{h \, X_{{\vartheta }_{\ell }}}(T)$ is also a Bohr-Sommerfeld torus 
$T^{\prime }_{{\mathbf{n}}^{\prime }}$, where ${\mathbf{n}}^{\prime }=(n_1, \ldots ,
n_{\ell -1}, n_{\ell }-1, n_{\ell +1}, \ldots , n_k)$.} \medskip

\noindent \textbf{Proof.} For simplicity we assume that $\ell =1$. Suppose
that the image of the curve $\gamma :[0,h] \rightarrow B: t \mapsto {\mathrm{%
e}}^{ \, X_{{\vartheta }_1}}(\rho (x_0)) $ lies in $V = {\psi }(U)$,
where $x_0 \in T = T_{\mathbf{n}}$. For $x \in T$ and $t \in [0,h]$ we have 
\begin{equation*}
X_{{\vartheta }_1} j_{\ell } = \left\{ 
\begin{array}{cccl}
X_{{\vartheta }_1} j_1 & \hspace{-5pt} = -\frac{\partial }{\partial j_1}j_1
& \hspace{-5pt}= -1, & \mbox{if $\ell =1$} \\ 
\rule{0pt}{12pt} X_{{\vartheta }_1} j_{\ell } & \hspace{-5pt} = -\frac{%
\partial }{\partial j_1}j_{\ell } & \hspace{-5pt} = \, \, \, \, 0, & %
\mbox{if $\ell =2, \ldots , k$}%
\end{array}
\right.
\end{equation*}
and $X_{{\vartheta }_1}{\vartheta }_{\ell } = -\frac{\partial }{\partial p_1}
{\vartheta }_{\ell } = 0$. Since $x \in T$ has action angle coordinates 
$(j_1(x), \ldots , j_k(x), {\vartheta }_1(x), \ldots , {\vartheta }_k(x))$ in 
$U$, the point ${\mathrm{e}}^{t\, {\vartheta }_1}(x)$ has action angle
coordinates $(j_1(x)-t, \ldots , j_k(x), {\vartheta }_1(x), \ldots , 
{\vartheta }_k(x))$. In particular, the point 
${\mathrm{e}}^{t\, X_{{\vartheta }_1}}(x_0)$ has action angle coordinates 
$(j_1(x_0)-t, \ldots ,j_k(x_0), {\vartheta }_1(x_0), \ldots , $ ${\vartheta }_k(x_0))$. So 
\begin{equation*}
({\mathrm{e}}^{h\, X_{{\vartheta }_1}})_{\ast }j_{\ell } = \left\{ 
\begin{array}{cl}
j_1 -h, & \mbox{if $\ell =1$} \\ 
j_{\ell }, & \mbox{if $\ell = 2, \ldots , k$}%
\end{array}
\right.
\end{equation*}
and $({\mathrm{e}}^{h\, X_{{\vartheta }_1}})_{\ast }{\vartheta }_{\ell } = {%
\vartheta }_{\ell }$ for $\ell =1, 2, \ldots , k$. Since $T$ is the
Bohr-Sommerfeld torus $T_{\mathbf{n}}$ we have $j_{\ell } = \int^1_0 j_{\ell
} \, \mathrm{d} {\vartheta }_{\ell } = n_{\ell }\, h$. Then 
\begin{align}
\int^1_0 ({\mathrm{e}}^{h\, X_{{\vartheta }_1}})_{\ast }j_1 \, \mathrm{d} %
\big( ({\mathrm{e}}^{h\, X_{{\vartheta }_1}})_{\ast }{\vartheta }_1 \big) &
= \int^1_0 (j_1-h) \mathrm{d} {\vartheta }_1  \notag \\
&\hspace{-1in} = j_1 -h = (n_1 -1) h.  \notag
\end{align}
So the torus ${\mathrm{e}}^{h \, X_{{\vartheta }_1}}(T)$ is a
Bohr-Sommerfeld torus $T_{{\mathbf{n}}^{\prime }}$ with ${\mathbf{n}}%
^{\prime }= (n_1, \ldots , n_{\ell -1},$ $n_{\ell }-1, n_{\ell +1}, \ldots ,
n_k)$. \medskip

Now consider the case when the image of the curve $\gamma :[0, h ]
\rightarrow B: t \mapsto {\mathrm{e}}^{t \, X_{{\vartheta }_1 }}(\rho (x_0))$
is not contained in $V$. This means that ${\mathrm{e}}^{t\, X_{{\vartheta }%
_1}}(U)$, where $U = {\rho }^{-1}(V)$, does not contain the torus $T$. Since 
${\mathrm{e}}^{t\, X_{{\vartheta }_1 }}$ is a $1$-parameter group of
symplectomorphisms of $(P, \omega )$, for every $t \in \mathbb{R} $ the
functions $\big( ({\mathrm{e}}^{t\, X_{{\vartheta }_1 }})_{\ast }j_{\ell }$,
with $\ell =1, \ldots ,k$ and $({\mathrm{e}}^{t\, X_{{\vartheta }_1
}})_{\ast }{\vartheta }_{\ell }, \, \ell =1, \ldots ,k $ are action angle
coordinates on $({\mathrm{e}}^{t\, X_{{\vartheta }_1 }})_{\ast }(U)$. Choose 
$\tau >0$ so that ${\mathrm{e}}^{\tau X_{{\vartheta }_1 }}(T) \subseteq U$.
Suppose that $h = \tau + \eta $, where $\eta \in [0, \tau )$. Observe that
for $t \in [0, \tau )$ the action angle coordinates $(j_1, \ldots , p_k, {%
\vartheta }_1,$ $\ldots , {\vartheta }_k)$ in $U$ satisfy 
\begin{equation*}
({\mathrm{e}}^{t X_{{\vartheta }_1 }})_{\ast }j_{\ell } = \left\{ 
\begin{array}{ll}
j_1 -t & \mbox{if $\ell =1$} \\ 
j_{\ell } & \mbox{if $\ell =2,3, \ldots ,k$}%
\end{array}
\right. \, \mathrm{and} \, \, \ ({\mathrm{e}}^{t \, X_{{\vartheta }_1
}})_{\ast }{\vartheta }_{\ell } = {\theta }_{\ell}.
\end{equation*}
Hence $({\mathrm{e}}^{\tau X_{{\vartheta }_1 }})_{\ast }j_1 = j_1 - \tau $
and 
\begin{align}
({\mathrm{e}}^{h \, X_{{\vartheta }_1 }})_{\ast }j_1 & = ({\mathrm{e}}%
^{(\tau + \eta )X_{{\vartheta }_1 }})_{\ast }j_1 = 
({\mathrm{e}}^{\tau \, X_{{\vartheta }_1 }})_{\ast }\, ({\mathrm{e}}^{\eta \, 
X_{{\vartheta }_1}})_{\ast }j_1  \notag \\
& = ({\mathrm{e}}^{\tau \, X_{{\vartheta }_1 }})_{\ast }(j_1 - \eta ) = 
({\mathrm{e}}^{\tau \, X_{{\vartheta }_1 }})_{\ast }(j_1) - \eta ,  \notag
\end{align}
because $\eta $ is constant. Moreover, 
\begin{align}
\int^1_0 \big( {\mathrm{e}}^{\tau \, X_{{\vartheta }_1} })_{\ast }j_1 \big) 
\mathrm{d} \big( {\mathrm{e}}^{\tau \, X_{{\vartheta }_1 }})_{\ast } {%
\vartheta }_1 \big) & = \int^1_0 (j_1 - \tau ) \, \mathrm{d} {\vartheta }_1
= \int^1_0 j_1 \, \mathrm{d} {\vartheta }_1 - \tau  \notag \\
& = j_1 - \tau .  \notag
\end{align}
Similarly, 
\begin{align}
\int^1_0 \big( {\mathrm{e}}^{h\, X_{{\vartheta }_1 }})_{\ast }j_1 \big) 
\mathrm{d} \big( {\mathrm{e}}^{h\, X_{{\vartheta }_1 }})_{\ast } {\vartheta }%
_1 \big) & = \int^1_0 \big( {\mathrm{e}}^{\tau \, X_{{\vartheta }_1
}})_{\ast }j_1 - \eta \big) \mathrm{d} \big( {\mathrm{e}}^{\tau \, X_{{%
\vartheta }_1 }})_{\ast } {\vartheta }_1 \big)  \notag \\
&\hspace{-1.25in} = \int^1_0 \big( {\mathrm{e}}^{\tau \, X_{{\vartheta }_1
}})_{\ast }j_1 \mathrm{d} \big( {\mathrm{e}}^{\tau \, X_{{\vartheta }_1
}})_{\ast } {\vartheta }_1 \big) - \eta \int^1_0 \mathrm{d} \big( {\mathrm{e}%
}^{\tau \, X_{{\vartheta }_1 }})_{\ast } {\vartheta }_1 \big)  \notag \\
&\hspace{-1.25in} = \int^1_0 p_1 \, \mathrm{d} {\vartheta }_1 - \tau - \eta
= \int^1_0 p_1 \, \mathrm{d} {\vartheta }_1 - h = (n_1 -1)h,  \notag
\end{align}
because $T$ is a Bohr-Sommerfeld torus $T_{\mathbf{n}}$ with quantum numbers 
$(n_1, \ldots , n_k)$. Thus ${\mathrm{e}}^{\, h \, X_{{\vartheta }_1}}(T)$
is a Bohr-Sommerfeld torus corresponding to the quantum numbers $(n_1-1,
n_2, \ldots , n_k)$. This argument may be extended to cover the case where $%
\hbar = k \tau + \eta $ for any positive integer $k$ and $\eta \in [0, \tau
) $. \hfill {\tiny $\blacksquare $}

\subsection{Monodromy}

%%%%%%%%%%%%%%

Suppose that $\mathcal{U} = {\{ U_i \} }_{i \in I}$ is a good covering of $P$
such that for every $i\in I$ the chart $(U_i, {\psi }_i)$ is the domain of a
local trivialization of the toral bundle $\rho : P \rightarrow B$,
associated to the fibrating toral polarization $D$ of $P$, given by the
local action angle coordinates 
\begin{equation*}
{\rho}_{\mid {U_i}}: U_i \rightarrow V_i \times {\mathbb{T}}^k: p \mapsto {\psi }_i(p) =
 (j^i, {\vartheta }^i) = (j^i_1, \ldots ,j^i_k, {\vartheta }^i_1,
\ldots , {\vartheta}^i_k)
\end{equation*}
with $( {\rho }_{\mid {U_i}} )_{\ast }({\omega }_{\mid {U_i}}) = 
\sum^k_{\ell =1}\mathrm{d} j^i_{\ell } \wedge \mathrm{d} {\vartheta }^i_{\ell }$. 
We suppose that
the set $\mathcal{S}$ of Bohr-Sommerfeld tori on $P$ has the local lattice
structure ${\{ {\mathbb{L}}_{U_i} \} }_{i \in I}$ of section 3.3. \medskip

\vspace{-.15in}Let $p$ and $p^{\prime }\in P$ and let $\gamma : [0,1]
\rightarrow P$ be a smooth curve joining $p$ to $p^{\prime }$. We can choose
a finite good subcovering ${\{ U_k \}}^N_{k=1}$ of $\mathcal{U}$ such that $%
\gamma ([0,1]) \subseteq \cup^N_{k=1}U_k$, where $\gamma (0) \subseteq U_1$
and $\gamma (1) \in U_N$. Using the fact that the local lattices ${\{ 
\mathbb{L}_{U_k} \} }^N_{k=1}$ are compatible, we can extend the local
action functions $j^1$ on $V_1 = {\psi}_1(U_1) \subseteq B$ to a local
action function $j^N$ on $V_N \subseteq B$. Thus using the connection $%
\mathcal{E}$, see corollary A.5, we may parallel transport a Bohr-Sommerfeld
torus $T_{\mathbf{n}} \subseteq U_1$ along the curve $\gamma $ to a
Bohr-Sommerfeld torus $T_{{\mathbf{n}}^{\prime }} \subseteq U_N$, see claim
3.4. The action function at $p^{\prime }$, in general depends on the path $%
\gamma $. If the holonomy group of the connection $\mathcal{E}$ on the
bundle $\rho : P \rightarrow B$ consists only of the identity element in $%
\mathrm{Gl}(k, \mathbb{Z} )$, then this extension process does not depend on
the path $\gamma $. Thus we have shown \medskip

\noindent \textbf{Claim 3.9} \textit{If $D$ is a fibrating toral
polarization of $(P, \omega )$ with fibration $\rho : P \rightarrow B$ and $B$ is simply connected, then there are global action angle coordinates on $P$ and the Bohr-Sommerfeld tori $T_{\mathbf{n}} \in \mathcal{%
S}$ have a unique quantum number $\mathbf{n}$. Thus the local lattice
structure of $\mathcal{S}$ is the lattice ${\mathbb{Z} }^k$.} \medskip

If the holonomy of the connection $\mathcal{E}$ on $P$ is not the identity
element, then the set $\mathcal{S}$ of Bohr-Sommerfeld tori is not a lattice
and it is not possible to assign a global labeling by quantum numbers to all
the tori in $\mathcal{S}$. This difficulty in assigning quantum numbers to
Bohr-Sommerfeld tori has been known to chemists since the early 1920s.
Modern papers illustrating it are \cite{winnewisser-et-al} and \cite%
{cushman-et-al}. We will give a concrete example where the connection $%
\mathcal{E}$ has nontrivial holonomy, namely, the spherical pendulum.
\medskip

\noindent \textbf{Example 3.10} The spherical pendulum is a completely
integrable Hamiltonian system $(H,J, T^{\ast }S^2, \mathrm{d} {\alpha }%
_{T^{\ast }S^2})$, where 
\begin{equation*}
T^{\ast }S^2 = \{ (q,p) \in T^{\ast }{\mathbb{R} }^3 \, 
\rule[-4pt]{.5pt}{13pt}\, \langle q, q \rangle = 1 \, \, \& \, \, \langle q,
p \rangle =0 \}
\end{equation*}
is the cotangent bundle of the $2$-sphere $S^2$ with $\langle \, \, , \, \,
\rangle $ the Euclidean inner product on ${\mathbb{R} }^3$. The Hamiltonian
is 
\begin{equation*}
H: T^{\ast }S^2 \rightarrow \mathbb{R} : (q,p) \mapsto \mbox{$\frac{%
\scriptstyle 1}{\scriptstyle 2}\,$} \langle p,p \rangle + \langle q, e_3
\rangle ,
\end{equation*}
where $e^T_3 = (0,0,1) \in {\mathbb{R} }^3$ and the $e_3$-component of
angular momentum is 
\begin{equation*}
J:T^{\ast }S^2 \rightarrow \mathbb{R} : (q,p) \mapsto q_1p_2 -q_2p_1.
\end{equation*}
The integral map of the spherical pendulum is 
\begin{equation*}
F: T^{\ast }S^2 \rightarrow \overline{R} \subseteq {\mathbb{R} }^2:(q,p)
\mapsto \big( H(q,p), J(q,p) \big) ,
\end{equation*}
see figure 1. Here $\overline{R}$ is the closure in ${\mathbb{R} }^2$ of the
set $R$ of regular values of the integral map $F$. The point $(1,0)$ is an
isolated critical value of $F$. So the set $R$ has the homotopy type of $S^1$
and is \emph{not} simply connected. Every fiber of $F_{\mid {F^{-1}(R)}}$ over a
point in $R$ is a $2$-torus, see \cite[chpt.V]{cushman-bates}. At every
point of $T^{\ast }S^2 \setminus F^{-1}(1,0)$ there are local action angle
coordinates. The fibers of $F$ corresponding to the dark points in figure $1$
are Bohr-Sommerfeld tori.

\begin{tabular}{l}
\setlength{\unitlength}{1pt} \\ 
\vspace{-.35in} \\ 
\hspace{1.1in}\includegraphics[width = 165pt]{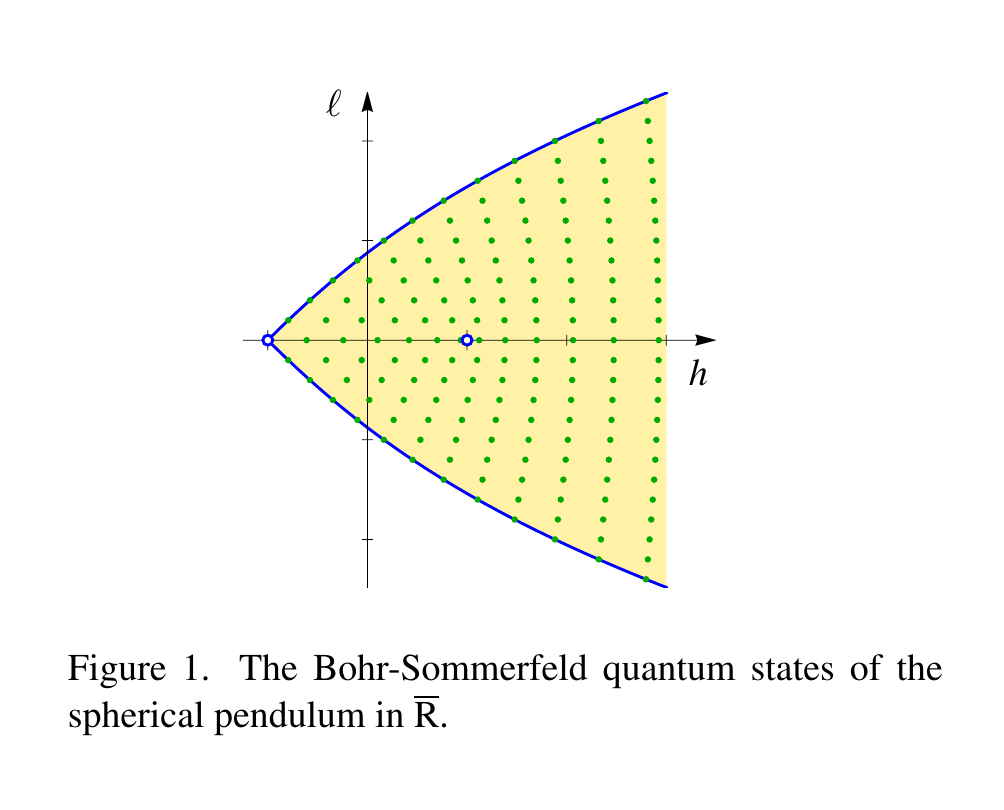} \\ 
\end{tabular}
\medskip

\vspace{-.15in}\noindent Since the are no global action angle coordinates,
the action function $j$ on $R$ is multi-valued. After encircling the point $%
(1,0)$, the quantum number of the torus represented by the upper right hand
vertices of the rectangle on the $h$-axis, see figure 2, becomes the quantum
number of the upper right hand vertex of the parallelogram formed by
applying {\tiny $%
\begin{pmatrix}
1 & 1 \\ 
0 & 1%
\end{pmatrix}%
$} to the original rectangle, which is the monodromy matrix $M$ of the
spherical pendulum.

\begin{tabular}{l}
\setlength{\unitlength}{1pt} \\ 
\vspace{-.3in} \\ 
\hspace{1in}\includegraphics[width = 165pt]{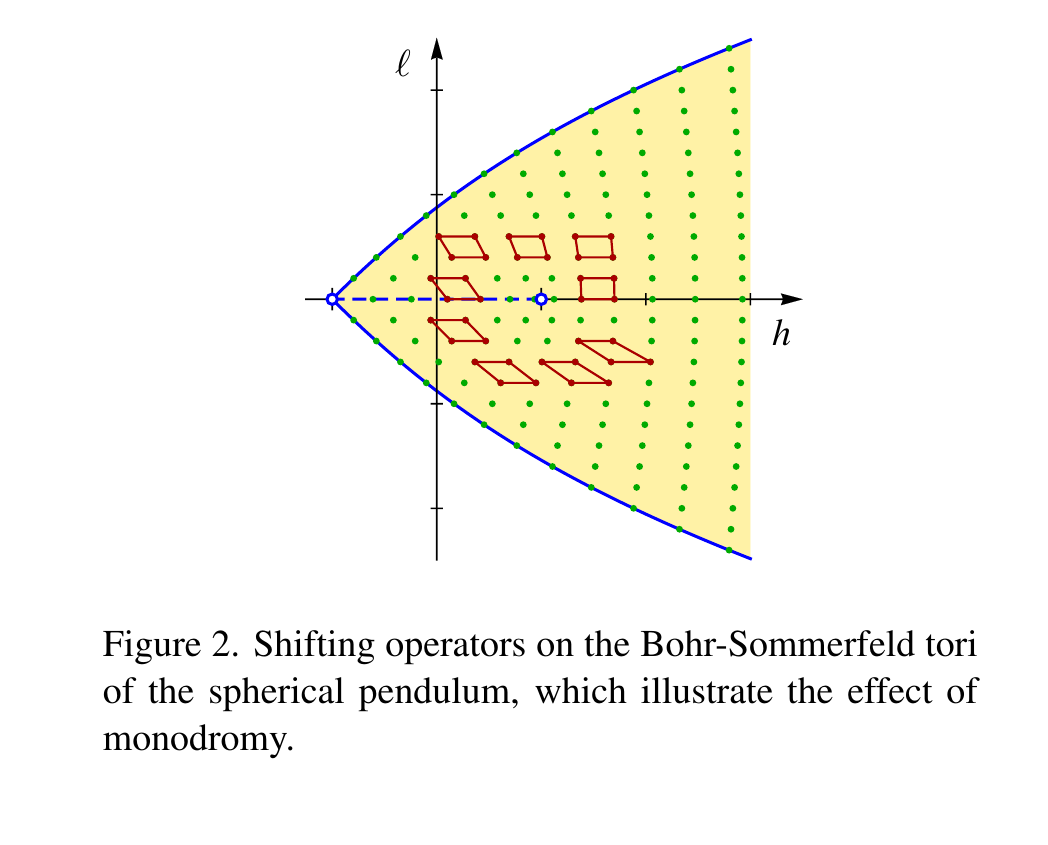} \\ 
\end{tabular}

The holonomy of the connection $\mathcal{E}$ is called the \emph{monodromy}
of the fibrating toral polarization $D$ on $(P, \omega )$ with fibration $%
\rho : P \rightarrow B$. \medskip

\noindent \textbf{Corollary 3.11} \textit{Let $\widetilde{B}$ be the
universal covering space of $B$ with covering map $\Pi : \widetilde{B}
\rightarrow B$. The monodromy map $M$, which is a nonidentity element
holonomy group of the connection $\mathcal{E}$ on the bundle $\rho $ sends
one sheet of the universal covering space to another sheet.} \medskip

\noindent \textbf{Proof.} Since the universal covering space $\widetilde{B}$
of $B$ is simply connected and we can pull back the symplectic manifold $(P,
\omega)$ and the fibrating toral distribution $D$ by the universal covering
map to a symplectic manifold $(\widetilde{P}, \widetilde{\omega })$ and a
fibrating toral distribution $\widetilde{D}$ with associated fibration $%
\widetilde{\rho}: \widetilde{P} \rightarrow \widetilde{B}$. The connection $%
\mathcal{E}$ on the bundle $\rho $ pulls back to a connection $\widetilde{%
\mathcal{E}}$ on the bundle $\widetilde{\rho}$. Let $\gamma $ be a closed
curve on $B$ and let $M$ be the holonomy of the connection $\mathcal{E}$ on $%
B$ along $\gamma $. Then $\gamma $ lifts to a curve $\widetilde{\gamma }$ on 
$\widetilde{B}$, which covers ${\gamma }$, that is, $\widetilde{\rho } \, %
\raisebox{2pt}{$\scriptstyle\circ \, $} \widetilde{\gamma } = \gamma $. Thus
parallel transport of a $k$-torus $T = {\mathbb{R} }^k/{\mathbb{Z} }^k$,
which is an integral manifold of the distribution $\widetilde{D}$, along the
curve $\widetilde{\gamma }$ gives a linear map $M$ of the lattice ${\mathbb{Z%
} }^k$ defining the $k$-torus $M(\widetilde{T})$. The map $M$ is the same as
the linear map $M$ of ${\mathbb{Z} }^k$ into itself given by parallel
transporting $T$, using the connection $\mathcal{E}$, along the closed $%
\gamma $ on $B$ because the connection $\widetilde{\mathcal{E}}$ is the pull
back of the connection $\mathcal{E}$ by the covering map $\rho $. The closed
curve $\gamma $ in $B$ represents an element of the fundamental group of $B$%
, which acts as a covering transformation on the universal covering space $%
\widetilde{B}$ that permutes the sheets (= fibers) of the universal covering
map $\widetilde{\Pi}$. \medskip

\noindent \textbf{Example 3.10 (continued)} In the spherical pendulum the
universal covering space $\widetilde{R}$ of $R \setminus \{ (1,0) \}$ is ${%
\mathbb{R} }^2$. If we cut $R$ by the line segment $\ell = \{ (h,0 ) \in R
\, \rule[-4pt]{.5pt}{13pt}\, h > 1 \}$, then $R^{\times } = R \setminus \ell 
$ is simply connected and hence represents one sheet of the universal
covering map of $R$. For more details on the universal covering map see \cite%
{cushman-sniatycki16}. The curve chosen in the example has holonomy $M=$%
{\tiny $%
\begin{pmatrix}
1 & 1 \\ 
0 & 1%
\end{pmatrix}%
$}. It gives a map of $\widetilde{R}$ into itself, which sends $R^{\times }$
to the adjacent sheet of the covering map. Thus we have a rule how the
labelling of the Bohr-Sommerfeld torus $T_{(n_1,n_2)}$, corresponding to $%
(h,j) \in R^{\times}$, changes when we go to an adjacent sheet, which covers 
$R^{\times}$, namely, we apply the matrix $M$ to the integer vector {\tiny $%
\begin{pmatrix}
n_1 \\ 
n_2%
\end{pmatrix}%
$}. Since our chosen curve generates the fundamental group of $R \setminus
\{ (1,0) \}$, we know what the quantum numbers of Bohr-Sommerfeld are for
any closed curve in $R \setminus \{ (1,0) \} $, which encircles the origin.
\hfill {\tiny $\blacksquare $}

\section{Appendix}

%%%%%%%%%%%%%

We return to study the symplectic geometry of a fibrating toral
polarization $D$ of the symplectic manifold $(P, \omega )$ in order to explain what
we mean by its local integral affine structure. \medskip

We assume that the integral manifolds ${\{ M_p \} }_{p \in P}$ of the
Lagrangian distribution $D$ on $P$ form a smooth manifold $B$ such that the
map 
\begin{equation*}
\rho : P \rightarrow B: p \mapsto M_p
\end{equation*}
is a proper surjective submersion. If the distribution $D$ has these
properties we refer to it as a \emph{fibrating polarization} of $(P, \omega
) $ with \emph{associated fibration} $\rho : P \rightarrow B$. \medskip

\noindent \textbf{Lemma A.1} \textit{Suppose that $D$ is a fibrating
polarization of $(P, \omega )$. Then the associated fibration $\rho : P
\rightarrow B$ has an Ehresmann connection $\mathcal{E}$ with parallel
translation. So the fibration $\rho : P \rightarrow B$ is locally trivial
bundle.} \medskip

\noindent \textbf{Proof} We construct the Ehresmann connection as follows.
For each $p \in P$ let $(U, \psi )$ be a \emph{Darboux chart} for $(P,
\omega )$. In other words, $({\psi }^{-1})^{\ast }({\omega }_{\mid U})$ is the
standard symplectic form ${\omega }_{2k}$ on $TV$, where $V = {\psi }(U)
\subseteq {\mathbb{R} }^{2k}$ with $\psi (p) = 0$. In more detail, for every 
$u \in U$ there is a frame $\varepsilon (u)$ of $P$ at $u$, whose image
under $T_u \psi $ is the frame $\varepsilon (v) = \big\{ \frac{\partial }{%
\partial x_1}\rule[-7pt]{.5pt}{15pt}\raisebox{-6pt}{$\, \scriptstyle{v}$},
\ldots , \frac{\partial }{\partial x_k}\rule[-7pt]{.5pt}{15pt}%
\raisebox{-6pt}{$\, \scriptstyle{v}$}, \, \frac{\partial }{\partial y_1}%
\rule[-7pt]{.5pt}{15pt}\raisebox{-6pt}{$\, \scriptstyle{v}$}, \ldots , \frac{%
\partial }{\partial y_k}\rule[-7pt]{.5pt}{15pt}%
\raisebox{-6pt}{$\,
\scriptstyle{v}$} \big\} $ of $T_vV = {\mathbb{R} }^{2k}$, where $v = \psi
(u)$, such that 
\begin{equation*}
{\omega }_{2k}(v) \big( \frac{\partial }{\partial x_i}%
\rule[-10pt]{.5pt}{24pt}\raisebox{-9pt}{$\,
\scriptstyle{v}$} , \frac{\partial }{\partial x_j}\rule[-10pt]{.5pt}{24pt}%
\raisebox{-9pt}{$\, \scriptstyle{v}$} \big) = {\omega }_{2k}(v) \big( \frac{%
\partial }{\partial y_i}\rule[-10pt]{.5pt}{24pt}%
\raisebox{-9pt}{$\,
\scriptstyle{v}$} , \frac{\partial }{\partial y_j}\rule[-10pt]{.5pt}{24pt}%
\raisebox{-9pt}{$\, \scriptstyle{v}$} \big) = 0
\end{equation*}
and 
\begin{equation*}
{\omega }_{2k}(v) \big( \frac{\partial }{\partial x_i}%
\rule[-10pt]{.5pt}{24pt}\raisebox{-9pt}{$\,
\scriptstyle{v}$} , \frac{\partial }{\partial y_j}\rule[-10pt]{.5pt}{24pt}%
\raisebox{-9pt}{$\, \scriptstyle{v}$} \big) = {\delta }_{ij}.
\end{equation*}

For $u \in M_p \cap U$, we see that ${\lambda }_v
= T_u\psi (T_uM_p)$ is a Lagrangian subspace of the symplectic vector space $%
\big( T_v V, {\omega }_{2k}(v) \big) $. Let ${\{ \frac{\partial }{\partial
z_j} \rule[-7pt]{.5pt}{15pt}\raisebox{-6pt}{$\, \scriptstyle{v}$} \} }%
^k_{j=1}$ be a basis of ${\lambda }_v$ with ${\{ \mathrm{d} z_j (v) \} }%
^k_{j=1}$ the corresponding dual basis of ${\lambda }^{\ast }_v$. Extend
each covector ${\mathrm{d} z}_j(v)$ by zero to a covector $\mathrm{d} Z_j(v)$
in $T^{\ast}_vV$, that is, extend the basis ${\{ \mathrm{d} z_j(v) \} }%
^k_{j=1}$ of ${\lambda }^{\ast }_v$ to a basis ${\{ \mathrm{d} Z_j(v) \} }%
^{k}_{j=1}$ of $T^{\ast }_v V$, where {\tiny $\left\{ 
\begin{array}{l}
\mathrm{d} Z_j(v)_{\mid {{\lambda }_v}} = \mathrm{d} z_j(v) , \, 
\mbox{for $j =1,
\ldots , k$} \\ 
\rule{0pt}{7pt} \mathrm{d} Z_j(v)_{\mid {{\lambda }_v}} = 0, \, 
\mbox{for $j =
k+1, \ldots , 2k$.}%
\end{array}
\right. $} Since ${\omega }^{\#}_{2k}(v): T_v V \rightarrow T^{\ast }_vV$ is
a linear isomorphism with inverse ${\omega }^{\flat}_{2k}(v)$ for every $v
\in V$, we see that the collection 
\begin{equation*}
{\{ \frac{\partial }{\partial w_j}\rule[-10pt]{.5pt}{21pt}%
\raisebox{-9pt}{$\, \scriptstyle{v}$} = {\omega }^{\flat}_{2k}(v) ( \mathrm{d%
} Z_j(v) ) \} }^{\raisebox{-4pt}{$\scriptstyle k$}}_ {\raisebox{4pt}{$%
\hspace{-2pt} \scriptstyle j=1$}}
\end{equation*}
of vectors in $T_vV$ spans an $k$-dimensional subspace ${\mu }_v$. We now
show that ${\mu }_v$ is a Lagrangian subspace of $\big( T_v V, {\omega }%
_{2k}(v) \big) $. By definition 
\begin{align}
{\omega }_{2k}(v) \big( \frac{\partial }{\partial w_i}%
\rule[-10pt]{.5pt}{24pt}\raisebox{-9pt}{$\, \scriptstyle{v}$}, \frac{%
\partial }{\partial w_j}\rule[-10pt]{.5pt}{24pt}%
\raisebox{-9pt}{$\,
\scriptstyle{v}$} \big) & = {\omega }^{\#}_{2k}(v) \big( \frac{\partial }{%
\partial w_i}\rule[-10pt]{.5pt}{24pt}\raisebox{-9pt}{$\, \scriptstyle{v}$} %
\big) \frac{\partial }{\partial w_j}\rule[-10pt]{.5pt}{24pt}%
\raisebox{-9pt}{$\, \scriptstyle{v}$} = \mathrm{d} Z_i(v) \frac{\partial }{%
\partial w_j}\rule[-10pt]{.5pt}{24pt}\raisebox{-9pt}{$\, \scriptstyle{v}$} =
0.  \notag
\end{align}
The last equality above follows because $\frac{\partial }{\partial w_j}%
\rule[-7pt]{.5pt}{15pt}\raisebox{-6pt}{$\, \scriptstyle{v}$} \notin {\lambda 
}_v$. To see this we note that 
\begin{align}
{\omega }_{2k}(v) \big( \frac{\partial }{\partial w_j}%
\rule[-10pt]{.5pt}{24pt}\raisebox{-9pt}{$\, \scriptstyle{v}$}, \frac{%
\partial }{\partial z_j}\rule[-10pt]{.5pt}{24pt}%
\raisebox{-9pt}{$\,
\scriptstyle{v}$} \big) & = \mathrm{d} Z_j(v) \frac{\partial }{\partial z_j}%
\rule[-10pt]{.5pt}{24pt}\raisebox{-9pt}{$\, \scriptstyle{v}$} = \mathrm{d}
z_j(v) \frac{\partial }{\partial z_j}\rule[-10pt]{.5pt}{24pt}%
\raisebox{-9pt}{$\, \scriptstyle{v}$} = 1.  \notag
\end{align}
The Lagrangian subspace ${\mu }_v$ is complementary to the Lagrangian
subspace ${\lambda }_v$, that is, $T_v V = {\lambda }_v \oplus {\mu }_v$ for
every $v \in V$. \medskip

Consequently, ${\mathrm{hor}}_u\, = T_v {\psi }^{-1} {\mu }_v$ is a
Lagrangian subspace of $\big( T_u U, \omega (u) \big) $, which is
complementary to the Lagrangian subspace $T_u M_p$. Since the mapping ${%
\mathrm{hor}}_{\mid U} : U \rightarrow TU: u \mapsto {\mathrm{hor}}_u$ is smooth
and has constant rank, it defines \linebreak a Lagrangian distribution ${%
\mathrm{hor}}_{\mid U}$ on $U$. Hence we have a Lagrangian distribution $\mathrm{%
hor}$ on $(P, \omega )$. Since $T_uM_p$ is the tangent space to the fiber ${%
\rho }^{-1}\big( \rho (p) \big) = M_p$, the distribution ${\mathrm{ver}}_{\mid U}:
U \rightarrow TU: u \mapsto {\mathrm{ver}}_u = T_uM_p = {\lambda }_v$
defines the vertical Lagrangian distribution $\mathrm{ver}$ on $P$. Because $%
{\mathrm{ver}}_u = \ker T_u \rho $, it follows that $T_u \rho ({\mathrm{hor}}%
_u) = T_{\rho (u)}B$. Hence the linear mapping $T_u\rho _{\mid {{\mathrm{hor}}}%
_u}: {\mathrm{hor}}_u \rightarrow T_{\rho (u)}B$ is an isomorphism. Since $%
T_p P = {\mathrm{hor}}_p \oplus {\mathrm{ver}}_p$ for every $p \in P$ and
the mapping $T_p\rho _{\mid {{\mathrm{hor}}_p}}: {\mathrm{hor}}_p \rightarrow
T_{\rho (p)}B$ is an isomorphism for every $p \in P$, the distributions $%
\mathrm{hor}$ and $\mathrm{ver}$ on $P$ define an \emph{Ehresmann connection}
$\mathcal{E}$ for the Lagrangian fibration $\rho : P \rightarrow B$. \hfill 
{\tiny $\blacksquare $} \medskip

Let $X$ be a smooth complete vector field on $B$ with flow ${\mathrm{e}}^{t
X}$. Because the linear mapping $T_p\rho _{\mid {{\mathrm{hor}}_p}}: {\mathrm{hor}}%
_p \rightarrow T_{\rho (p)}B$ is bijective, there is a unique smooth vector
field $\mathrm{lift} X$ on $P$, called the \emph{horizontal lift} of $X$,
which is $\rho$-related to $X$, that is, $T_p\rho \, \mathrm{lift} X(p) = X%
\big( \rho (p) \big) $ for every $p \in P$. Let ${\mathrm{e}}^{t\, \mathrm{%
lift} X}$ be the flow of $\mathrm{lift} X$. Then $\rho ({\mathrm{e} }^{t \, 
\mathrm{lift} X}) = {\mathrm{e}}^{t X} (\rho (p))$. Let $\sigma : W
\subseteq B \rightarrow P$ be a smooth local section of the bundle $\rho : P
\rightarrow B$. Define the \emph{covariant derivative} ${\nabla }_X\sigma $
of $\sigma $ with respect to the vector field $X$ by 
\begin{equation*}
({\nabla }_X \sigma )(w) = 
\mbox{${\displaystyle \frac{\dee }{\dee t}}
\rule[-10pt]{.5pt}{25pt} \raisebox{-10pt}{$\, {\scriptstyle t=0}$}$} {%
\mathrm{e}}^{-t\, \mathrm{lift}\, X}\big( \sigma ({\mathrm{e}}^{t \, X} (w)) %
\big)
\end{equation*}
for all $w \in W$. Because the bundle projection map $\rho $ is proper, 
\emph{parallel transport} of each fiber of the bundle $\rho : P \rightarrow
B $ by the flow of $\mathrm{lift}X$ is defined as long as the flow of $X$ is
defined. Because the Ehresmann connection $\mathcal{E}$ has parallel
transport, the bundle presented by $\rho $ is locally trivial, see \cite[%
p.378--379]{cushman-bates}.\medskip

\noindent \textbf{Claim A.2} \textit{If $D$ is a fibrating polarization of
the symplectic manifold $(P, \omega )$, then for every $p \in P$ the
integral manifold of $D$ through $p$ is a smooth \linebreak Lagrangian
submanifold of $P$, which is an $k$-torus $T$. In fact $T$ is the fiber over 
$\rho (p)$ of the associated fibration $\rho : P \rightarrow B$.} \medskip

We say that $D$ is a \emph{fibrating toral polarization} of $(P, \omega )$
if it satisfies the hypotheses of claim A.2. The proof of claim A.2 requires
several preparatory arguments. \medskip

Let $f \in C^{\infty}(B)$. Then ${\rho }^{\ast }f \in C^{\infty}(P)$. Let $%
X_{{\rho }^{\ast }f}$ be the Hamiltonian vector field on $(P, \omega )$ with
Hamiltonian ${\rho }^{\ast }f$. We have \medskip

\noindent \textbf{Lemma A.3} \textit{Every fiber of the locally trivial
bundle $\rho : P \rightarrow B$ is an invariant manifold of the Hamiltonian
vector field $X_{{\rho }^{\ast }f}$.} \medskip

\noindent \textbf{Proof.} We need only show that for every $p \in P$ and
every $q \in M_p$, we have $X_{{\rho }^{\ast }f}(q) \in T_qM_p$. Let $Y$ be
a smooth vector field on the integral manifold $M_p$ with flow ${\mathrm{e}}%
^{t Y}$. Then 
\begin{equation*}
{\rho }^{\ast }f\big( {\mathrm{e}}^{t Y}(q) \big) = f\big( \rho ({\mathrm{e}}%
^{t Y}(q)) \big) = f\big( \rho (p) \big) ,
\end{equation*}
since ${\mathrm{e}}^{t Y}$ maps $M_p$ into itself. So 
\begin{align}
0 & = 
\mbox{${\displaystyle \frac{\dee }{\dee t}}
\rule[-10pt]{.5pt}{25pt} \raisebox{-10pt}{$\, {\scriptstyle t=0}$}$} {\rho }%
^{\ast }f\big( {\mathrm{e}}^{t Y}(q) \big) = L_Y ({\rho }^{\ast }f)(q) = 
\mathrm{d} \big( {\rho }^{\ast }f \big) (q) Y(q)  \notag \\
& = -{\omega }(q) \big( X_{{\rho }^{\ast }f}(q), Y(q) \big).  \notag
\end{align}
But $T_qM_p$ is a Lagrangian subspace of the symplectic vector space $(T_qP, 
{\omega }(q) )$. Consequently, $X_{{\rho }^{\ast }f}(q) \in T_q M_p$. \quad %
\mbox{\tiny $\blacksquare $} \medskip

Since the mapping $\rho : P \rightarrow B$ is surjective and proper, for
every $b \in B$ the fiber ${\rho }^{-1}(b)$ is a smooth compact submanifold
of $P$. Hence the flow ${\mathrm{e}}^{t \, X_{{\rho }^{\ast }f}}$ of the
vector field $X_{{\rho }^{\ast} f}$ is defined for all $t \in \mathbb{R} $.
\medskip

\noindent \textbf{Lemma A.4} \textit{Let $f$, $g \in C^{\infty}(B)$. Then $%
\{ {\rho }^{\ast }f , {\rho }^{\ast }g \} =0$.} \medskip

\noindent \textbf{Proof} For every $p \in P$ and every $q \in M_p$ from
lemma A.3 it follows that $X_{{\rho }^{\ast }f}(q)$ and $X_{{\rho }^{\ast
}g}(q)$ lie in $T_qM_p$. Because $M_p$ is a Lagrangian submanifold of $(P,
\omega )$, we get 
\begin{equation}
0 = \omega (q)\big( X_{{\rho }^{\ast }g}(q), X_{{\rho }^{\ast }f}(q) \big) =
\{ {\rho }^{\ast }f, {\rho }^{\ast }g \} (q).  \label{eq-s3ss5newone}
\end{equation}
Since $P = \amalg_{p \in P} M_p$, we see that (\ref{eq-s3ss5newone}) holds
for every $p \in P$. \quad \mbox{\tiny $\blacksquare $} \medskip

\noindent \textbf{Proof of claim A.2} From lemma A.4 it follows that $\big( {%
\rho }^{\ast }( C^{\infty}(B)), \{ \, \, , \, \, \}, \cdot \big)$ is an
abelian subalgebra $\mathfrak{t}$ of the Poisson algebra $(C^{\infty}(P), \{
\, \, , \, \, \} , \cdot )$. Since the bundle projection mapping $\rho : P
\rightarrow B$ is surjective and $\dim B =k$, the algebra $\mathfrak{t}$ has 
$k$ generators, say, ${\ \{ {\rho }^{\ast }f_i \} }^k_{i=1}$, whose
differentials at $q$ span $T_q({\rho }^{-1}(b))$ for every $b \in B$ and
every $q \in {\rho }^{-1}(b)$. Using the flow ${\mathrm{e}}^{t\, X_{{\rho }%
^{\ast }f_i}}$ of the Hamiltonian vector field $X_{{\rho }^{\ast}f_i}$ on $%
(P, \omega )$ define the ${\mathbb{R} }^k$-action 
\begin{equation}
\Phi : {\mathbb{R} }^k \times P \rightarrow P; \big( \mathbf{t} = (t_1,
\ldots , t_k), p \big) \mapsto \big( {\mathrm{e}}^{t_1 X_{{\rho }^{\ast
}f_1}}(p), \ldots , {\mathrm{e}}^{t_k X_{{\rho }^{\ast }f_k}}(p) \big)
\label{eq-s3ss5newtwo}
\end{equation}
Since ${\mathop{\rm span}\nolimits }_{1 \le i \le k}\{ X_{{\rho }^{\ast
}f_i}(q) \} = T_q({\rho}^{-1}(b))$ and each fiber is connected, being an
integral manifold of the distribution $D$, it follows that the ${\mathbb{R} }%
^k$-action $\Phi $ is transitive on each fiber ${\rho }^{-1}(b)$ of the
bundle $\rho : P \rightarrow B$. Thus ${\rho }^{-1}(b)$ is diffeomorphic to $%
{\mathbb{R} }^k/P_q$, where $P_q = \{ \mathbf{t} \in {\mathbb{R} }^k \, 
\rule[-4pt]{.5pt}{13pt}\, \, {\Phi }_{\mathbf{t}}(q) = q \} $ is the
isotropy group at $q$. If $P_q = \{ 0 \} $ for some $q \in P$, then the
fiber ${\rho }^{-1}\big( \rho (q) \big) $ would be diffeomorphic to ${%
\mathbb{R} }^k/P_q = {\mathbb{R} }^k$. But this contradicts the fact the
every fiber of the bundle $\rho : P \rightarrow B$ is compact. Hence $P_q
\ne \{ 0 \} $ for every $q \in P$. Since ${\mathbb{R} }^k/P_q$ is
diffeomorphic to ${\rho }^{-1}(b) $, they have the same dimension, namely, $%
k $. Hence $P_q$ is a zero dimensional Lie subgroup of ${\mathbb{R} }^k$.
Thus $P_q$ is a rank $k$ lattice ${\mathbb{Z} }^k$. So the fiber ${\rho }%
^{-1}(b)$ is ${\mathbb{R} }^k / {\mathbb{Z} }^k$, which is an \emph{affine} $%
k$-torus ${\mathbb{T} }^k$. \quad \mbox{\tiny $\blacksquare $} \medskip

We now apply the action angle theorem \cite[chpt.IX]{cushman-bates} to the
fibrating toral Lagrangian polarization $D$ of the symplectic manifold $(P,
\omega )$ with \linebreak associated toral bundle $\rho : P \rightarrow B$
to obtain a more precise description of the Ehresmann connection $\mathcal{E}
$ constructed in lemma A.2. For every $p\in P$ there is an open neighborhood 
$U$ of the fiber ${\rho }^{-1}\big( \rho (p) \big)$ in $P$ and a symplectic
diffeomorphism 
\begin{equation*}
\begin{array}{l}
\psi : U = {\rho }^{-1}(V) \subseteq P \rightarrow V \times {\mathbb{T} }^k
\subseteq {\mathbb{R} }^k \times {\mathbb{T} }^k: \\ 
\hspace{.5in} u \mapsto (j, \vartheta ) = (j_1, \ldots , j_k, {\vartheta }_1, \ldots , 
{\vartheta }_k)%
\end{array}%
\end{equation*}
such that 
\begin{equation*}
\rho _{\mid U} : U \subseteq P \rightarrow V \subseteq {\mathbb{R} }^k:u \mapsto ({%
\pi }_1 \, \raisebox{2pt}{$\scriptstyle\circ \, $} \psi )(u) = j ,
\end{equation*}
is the momentum mapping of the Hamiltonian ${\mathbb{T}}^k$-action on $(U,
\omega _{\mid U})$. Here ${\pi }_1: V \times {\mathbb{T} }^k \rightarrow V:(j,
\vartheta ) \rightarrow j$. Thus the bundle $\rho : P \rightarrow B$ is locally a
principal ${\mathbb{T} }^k$-bundle. Moreover, we have 
$({\psi }^{-1})^{\ast} \omega _{\mid U} = \sum^k_{i=1} \mathrm{d} j_i \wedge 
\mathrm{d} {\vartheta }_i$.
\medskip

\noindent \textbf{Corollary A.5} \textit{Using the chart $(U, \psi )$ for
action angle coordinates $(j, \phi )$, the Ehresmann connection 
${\mathcal{E}}_{\mid U}$ gives an Ehresmann connection 
${\mathcal{E}}_{\mid {V\times {\mathbb{T}}^n}}$
on the bundle ${\pi }_1: V \times {\mathbb{T} }^k \rightarrow V$ defined by} 
\begin{equation*}
{\mathrm{ver}}_v = {\mathop{\rm span}\nolimits }_{1\le i \le k}\{ \frac{%
\partial }{\partial {\vartheta }_i}\rule[-9pt]{.5pt}{18pt} 
\raisebox{-8pt}{$\,
\scriptstyle v = \psi (u)$} \} \, \, \, \mathrm{and} \, \, \, {\mathrm{hor}}%
_v = {\mathop{\rm span}\nolimits }_{1\le i \le k}\{ \frac{\partial }{\partial j_i}
\rule[-9pt]{.5pt}{18pt} \raisebox{-8pt}{$\, \scriptstyle v =\psi (u)$} \} .
\end{equation*}
\medskip

\noindent \textbf{Proof} This follows because $T_u\psi \big( {\mathrm{ver}}%
_u \big) = {\mathop{\rm span}\nolimits }_{1\le i \le k}\{ \frac{\partial }{%
\partial {\vartheta }_i}\rule[-9pt]{.5pt}{18pt} 
\raisebox{-8pt}{$\, \scriptstyle
v = \psi (u)$} \} $ and \linebreak $T_p \psi \big( {\mathrm{hor}}_u \big) = {%
\mathop{\rm span}\nolimits }_{1\le i \le k}\{ \frac{\partial }{\partial j_i}%
\rule[-9pt]{.5pt}{18pt} \raisebox{-8pt}{$\, \scriptstyle v = \psi (u)$} \} $
for every $u \in U$. From the preceding equations for every $u \in U$ we
have ${\mathrm{ver}}_u = {\mathop{\rm span}\nolimits }_{1\le i \le k}\{ X_{{%
\rho }^{\ast }(j_i)}(u) \} $ and ${\mathrm{hor}}_u = {\mathop{\rm span}%
\nolimits }_{1 \le i \le k}\{ X_{ ( {\pi }_2 \raisebox{-2pt}{$\comp$} \psi)^{\ast }
(-{\vartheta }_i )} (u) \} $. Here $\pi _2: V \times {\mathbb{T} }^k
\rightarrow {\mathbb{T} }^k:( j, \varphi ) \mapsto \varphi $. \quad 
\mbox{\tiny
$\blacksquare $} \medskip

\noindent \textbf{Corollary A.6 }\textit{The Ehresmann connection $\mathcal{E%
}$ on the locally trivial toral Lagrangian bundle $\rho : P \rightarrow B$
is flat, that is, ${\nabla }_X \sigma =0$ for every smooth vector field $X$
on $B$ and every local section $\sigma $ of $\rho : P \rightarrow B$.}
\medskip

\noindent \textbf{Proof} In action angle coordinates a local section section 
$\sigma $ of the bundle $\rho : P \rightarrow B$ is given by $\sigma : V
\rightarrow V \times {\mathbb{T} }^k: j \mapsto \big( j, \sigma (j) \big) $.
Let $X = \frac{\partial }{\partial j_{\ell }}$ for some $1 \le \ell \le k$
with flow ${\mathrm{e}}^{t\, X} $. Let $\mathrm{lift}X$ be the horizontal
lift of $X$ with respect to the Ehresmann connection ${\mathcal{E}}_{V
\times {\mathbb{T} }^k}$ on the bundle ${\pi }_1: V \times {\mathbb{T} }^k
\rightarrow V$. So for every $j \in V$ we have 
\begin{align}
( {\nabla }_X \sigma )(j) & = 
\mbox{${\displaystyle \frac{\dee }{\dee t}}
\rule[-10pt]{.5pt}{25pt} \raisebox{-10pt}{$\, {\scriptstyle t=0}$}$} {%
\mathrm{e}}^{t \, \mathrm{lift} X} \big( \sigma ({\mathrm{e}}^{-t X} (j)) %
\big)  \notag \\
& = 
\mbox{${\displaystyle \frac{\dee }{\dee t}}
\rule[-10pt]{.5pt}{25pt} \raisebox{-10pt}{$\, {\scriptstyle t=0}$}$} {%
\mathrm{e} }^{t \, \mathrm{lift} X} \big( \sigma ( j(-t) ) \big) , \quad %
\mbox{where ${\mathrm{e}}^{t X}(j) = j(t)$}  \notag \\
& = 
\mbox{${\displaystyle \frac{\dee }{\dee t}}
\rule[-10pt]{.5pt}{25pt} \raisebox{-10pt}{$\, {\scriptstyle t=0}$}$} {%
\mathrm{e}}^{t\, \mathrm{lift} X} \big( j , \sigma (j) \big), \quad %
\mbox{since $j_i$ for $1 \le i \le n$ are integrals of $X$}  \notag \\
& = 
\mbox{${\displaystyle \frac{\dee }{\dee t}}
\rule[-10pt]{.5pt}{25pt} \raisebox{-10pt}{$\, {\scriptstyle t=0}$}$} \big( %
j(t) , \sigma (j(t)) \big), \quad 
\mbox{since 
${\pi }_1 \big( {\mathrm{e}}^{t\, \mathrm{lift} X} ( j , \sigma (j) ) \big) = 
{\mathrm{e}}^{t X}(j)$}  \notag \\
& = 0.  \notag
\end{align}
This proves the corollary, since every vector field $X$ on $W \subseteq B$
may be written as $\sum^k_{i=1} c_i(j) \frac{\partial }{\partial j_i}$ for
some $c_i \in C^{\infty}(W)$ and the flow ${\{ {\varphi }^{\, j_i}_t \} }%
^k_{i=1}$ of ${\{ \frac{\partial }{\partial j_i} \} }^k_{i=1}$ on $V$
pairwise commute. \quad \mbox{\tiny $\blacksquare $} \bigskip

\noindent \textbf{Claim A.7 }\textit{Let $\rho : P \rightarrow B$ be a
locally trivial toral Lagrangian bundle, where $(P, \omega )$ is a smooth
symplectic manifold. Then the smooth manifold $B$ has an \emph{integral
affine} structure. In other words, there is a good open covering ${\{ W_i \} 
}_{i \in I}$ of $B$ such that the overlap maps of the coordinate charts $%
(W_i , {\varphi }_i)$ given by 
\begin{equation*}
{\varphi }_{i\ell } = {\varphi }_{\ell } \, \raisebox{2pt}{$\scriptstyle%
\circ \, $} {\varphi }^{-1}_i: V_i \cap V_{\ell } \subseteq {\mathbb{R} }^k
\rightarrow V_i \cap V_{\ell } \subseteq {\mathbb{R} }^k,
\end{equation*}
where ${\varphi }_i(W_i) = V_i$, have derivative $D{\varphi }_{i \ell }(v)
\in \mathrm{Gl}(k, \mathbb{Z} )$, which does not \linebreak depend on $v \in
V_i \cap V_{\ell }$.} \medskip

\noindent \textbf{Proof} Cover $P$ by $\mathcal{U} = {\ \{ U_i \} }_{i \in
I} $, where $(U_i , {\psi }_i)$ is an action angle coordinate chart. Since
every open covering of $P$ has a good refinement, we may assume that $%
\mathcal{U}$ is a good covering. Let $W_i = \rho (U_i)$. Then $\mathcal{W} = 
{\ \{ W_i \} }_{i \in I}$ is a good open covering of $B$ and $(W_i, {\varphi 
}_i = {\pi }_1 \, \raisebox{2pt}{$\scriptstyle\circ \, $} {\psi }_i)$ is a
coordinate chart for $B$. By construction of action angle coordinates, in $%
V_i \cap V_{\ell }$ the overlap map ${\varphi }_{i \ell }$ sends the action
coordinates $j^i$ in $V_i \cap V_{\ell }$ to the action coordinates $j^{\ell
}$ in $V_i \cap V_{\ell }$. The period lattices $P_{{\psi }^{-1}_i(j^i)}$
and $P_{{\psi }^{-1}_{\ell }(j^{\ell })}$ are equal since for some $p \in
W_i \cap W_{\ell}$ we have ${\psi }_i(p) = j^i$ and ${\psi }_{\ell }(p) =
j^{\ell }$. Moreover, these lattices do not depend on the point $p$. Thus
the derivative $D{\varphi }_{i\ell}(j)$ sends the lattice ${\mathbb{Z} }^k$
spanned by ${\ \{ \frac{\partial }{\partial j^i} \rule[-8pt]{.5pt}{18pt}%
\raisebox{-6pt}{$\, \scriptstyle j$} \} }^k_{i=1}$ into itself. Hence for
every $j \in W_i \cap W_{\ell }$ the matrix of $D{\varphi }_{i\ell}(j)$ has
integer entries, that is, it lies in $\mathrm{Gl}(k, \mathbb{Z})$ and the
map $j \mapsto D{\varphi }_{i\ell}(j)$ is continuous. But $\mathrm{Gl}(k, 
\mathbb{Z} )$ is a discrete subgroup of the Lie group $\mathrm{Gl}(k, 
\mathbb{R} )$ and $W_i \cap W_{\ell }$ is connected, since $\mathcal{W}$ is
a good covering. So $D{\varphi }_{i\ell}(j)$ does not depend on $j \in W_i
\cap W_{\ell }$. \quad \mbox{\tiny $\blacksquare$} \medskip

\noindent \textbf{Corollary A.8} \textit{Let $\gamma :[0,1] \rightarrow B$
be a smooth closed curve in $B$. Let $P_{\gamma }: [0,1] \rightarrow P$ be
parallel translation along $\gamma $ using the Ehresmann connection $%
\mathcal{E}$ on the bundle $\rho : P \rightarrow B$. Then the holonomy group
of the $k$-toral fiber $T_{\gamma (0)} = {\mathbb{T} }^k$ is induced by the
group $\mathrm{Gl}(k , \mathbb{Z} ) \ltimes {\mathbb{Z} }^k$ of affine $%
\mathbb{Z} $-linear maps of ${\mathbb{Z} }^k$ into itself.}

\end{document}